\documentclass[10pt]{amsart}
\usepackage{amsmath}
\usepackage{amsthm}
\usepackage[all]{xy}

\theoremstyle{plain}
\newtheorem{thm}{Theorem}
\newtheorem{prop}[thm]{Proposition}
\newtheorem{lem}[thm]{Lemma}

\newtheorem{conj}[thm]{Conjecture}
\newtheorem{ithm}{Theorem}
\newtheorem{iconj}[ithm]{Conjecture}

\DeclareMathOperator{\aut}{Aut}
\DeclareMathOperator{\depth}{depth}
\DeclareMathOperator{\coker}{coker}
\DeclareMathOperator{\der}{Der}
\DeclareMathOperator{\enm}{End}
\DeclareMathOperator{\ext}{Ext}
\DeclareMathOperator{\glg}{GL}

\DeclareMathOperator{\hmm}{Hom}

\DeclareMathOperator{\im}{im}

\DeclareMathOperator{\rk}{rk}
\DeclareMathOperator{\slg}{SL}
\DeclareMathOperator{\spec}{Spec}
\DeclareMathOperator{\syz}{syz}

\newcommand{\diff}{\mathsf{Diff}}
\newcommand{\g}{\mathsf{g}}
\newcommand{\Gm}{\Gamma}
\newcommand{\kdf}[1]{\Omega_{#1}}
\newcommand{\loc}{S^{-1}}
\newcommand{\mc}{\mathsf{MC}}
\newcommand{\omc}{\mathsf{\Omega MC}}
\newcommand{\mic}{\mathsf{MIC}}
\newcommand{\omic}{\mathsf{\Omega MIC}}

\newcommand{\n}{\mathbf{N_0}}

\newcommand{\z}{\mathbf{Z}}
\newcommand{\q}{\mathbf{Q}}

\begin{document}

\title[Connections on modules over singularities]{Connections on modules
over singularities of finite CM representation type}
\author{Eivind Eriksen}
\address{Oslo University College}
\email{eeriksen@hio.no}
\author{Trond St{\o}len Gustavsen}
\address{University of Oslo}
\email{stolen@math.uio.no}
\date{\today}

\maketitle

\begin{abstract}
Let $A$ be a commutative $k$-algebra, where $k$ is an algebraically
closed field of characteristic $0$, and let $M$ be an $A$-module. We
consider the following question: Under what conditions on $A$ and $M$
is it possible to find a connection $\nabla: \der_k(A) \to \enm_k(M)$
on $M$?

We consider maximal Cohen-Macaulay (MCM) modules over complete CM
algebras that are isolated singularities, and usually assume that the
singularities have finite CM representation type. It is known that any
MCM module over a simple singularity of dimension $d \le 2$ admits an
integrable connection. We prove that an MCM module over a simple
singularity of dimension $d \ge 3$ admits a connection if and only if
it is free. Among singularities of finite CM representation type, we
find examples of curves with MCM modules that do not admit connections,
and threefolds with non-free MCM modules that admit connections.

Let $A$ be a singularity not necessarily of finite CM representation
type, and consider the condition that $A$ is a Gorenstein curve or a
$\q$-Gorenstein singularity of dimension $d \ge 2$. We show that
this condition is sufficient for the canonical module $\omega_A$ to
admit an integrable connection, and conjecture that it is also
necessary. In support of the conjecture, we show that if $A$ is a
monomial curve singularity, then the canonical module $\omega_A$
admits an integrable connection if and only if $A$ is Gorenstein.
\end{abstract}

\section*{Introduction}

Let $k$ be an algebraically closed field of characteristic $0$, and
let $A$ be a commutative $k$-algebra. For any $A$-module $M$, we
consider the notion of a \emph{connection} on $M$, i.e. an
$A$-linear homomorphism
    \[ \nabla: \der_k(A) \to \enm_k(M) \]
such that $\nabla_D(am) = a \nabla_D(m) + D(a) m$ for all $D \in
\der_k(A)$, $a \in A$ and $m \in M$. A connection is integrable if
it is a Lie algebra homomorphism. The present paper is devoted to
the following question: Under what conditions on $A$ and $M$ is it
possible to find a connection on $M$?

We consider a complete local CM $k$-algebra $A$ with residue field $k$
that is an isolated singularity, and a maximal Cohen-Macaulay (MCM)
$A$-module $M$. Moreover, we usually assume that $A$ has finite CM
representation type, i.e. the number of isomorphism classes of
indecomposable MCM $A$-modules is finite.

If $A$ is a hypersurface, then $A$ has finite CM representation type if
and only if it is a simple singularity. By convention, $A =
k[[x]]/(x^{n+1})$ for $n \ge 1$ are the simple singularities (of type
$A_n$) of dimension zero. It is known that if $A$ is a simple
singularity of dimension $d \le 2$, then any MCM $A$-module admits an
integrable connection, see section \ref{s:mcm-conn}. We prove the
following result:

\begin{ithm}
Let $A$ be the complete local ring of a simple singularity of
dimension $d \ge 3$. Then an MCM $A$-module $M$ admits a connection
if and only if $M$ is free.
\end{ithm}

A Gorenstein singularity with finite CM representation type is a
hypersurface, and there is a complete classification of non-Gorenstein
singularities of finite CM representation type in dimension $d \le 2$:
\begin{enumerate}
\item The curve singularities $D^s_n$ for $n \ge 2$ and $E^s_6, E^s_7,
E^s_8$.
\item The quotient surface singularities that are non-Gorenstein.
\end{enumerate}
On the other hand, the classification is not complete in higher
dimensions. The only known examples of non-Gorenstein singularities of
finite CM representation type in dimension $d \ge 3$ are the following
threefolds:
\begin{enumerate}
\setcounter{enumi}{2}
\item The quotient threefold singularity of type $\frac{1}{2} (1,1,1)$.
\item The threefold scroll of type $(2,1)$.
\end{enumerate}
Definitions of the singularities in (1)-(4) are given in section
\ref{s:mcm-conn}.

Among the curve singularities of finite CM representation type, we
find examples of singularities with MCM modules that do not admit
connections. In fact, the canonical module $\omega_A$ does not admit
a connection when $A$ is the complete local ring of one of the
singularities $E^s_6, E^s_7, E^s_8$. Moreover, it seems that the
same holds for the singularities $D^s_n$ for all $n \ge 2$.

However, for the surface singularities of finite CM representation
type, all MCM modules admit connections. This is a consequence of
the fact that these singularities are quotient singularities, and
that all MCM modules are induced by group representations. It is
perhaps more surprising that it is difficult to find examples of MCM
modules over surface singularities that do not admit connections,
even when the singularities have infinite CM representation type.

In dimension $d \ge 3$, we would like to find examples of non-free
MCM modules that admit connections. We find that over the threefold
scroll of type $(2,1)$, no non-free MCM modules admit connections.
However, over the threefold quotient singularity of type
$\frac{1}{2} (1,1,1)$, the canonical module admits an integrable
connection.

Let $A$ be a singularity not necessarily of finite CM representation
type. Let us consider the condition that $A$ is a Gorenstein curve or a
$\q$-Gorenstein singularity of dimension $d \ge 2$. We show that this
condition is sufficient for the canonical module $\omega_A$ to admit an
integrable connection. In accordance with all our results for
singularities of finite CM representation type, we make the following
conjecture:

\begin{iconj}
Let $A$ be the complete local ring of a singularity of dimension $d \ge
1$. Then the canonical module $\omega_A$ admits a connection if and
only if $A$ is a Gorenstein curve or a $\q$-Gorenstein singularity of
dimension $d \ge 2$.
\end{iconj}

Let us consider a monomial curve singularity $A$ not necessarily of
finite CM representation type. We show that if $A$ is Gorenstein, then
any gradable rank one MCM $A$-module admits an integrable connection.
We also prove the following theorem, in support of our conjecture:

\begin{ithm}
Let $A$ be the complete local ring of a monomial curve singularity.
Then the canonical $A$-module $\omega_A$ admits a connection if and
only if $A$ is Gorenstein.
\end{ithm}

\section{Basic definitions}
\label{s:obstr}

Let $k$ be an algebraically closed field of characteristic $0$, and let
$A$ be a commutative $k$-algebra. A \emph{Lie-Rinehart algebra} of
$A/k$ is a pair $(\g, \tau)$, where $\g$ is an $A$-module and a $k$-Lie
algebra, and $\tau: \g \to \der_k(A)$ is a morphism of $A$-modules and
$k$-Lie algebras, such that
    \[ [D, a D'] = a [D,D'] + \tau_D(a) \; D' \]
for all $D,D' \in \g$ and all $a \in A$, see Rinehart \cite{ri63}. A
Lie-Rinehart algebra is the algebraic analogue of a \emph{Lie
algebroid}, and it is also known as a Lie pseudo-algebra or a
Lie-Cartan pair.

When $\g$ is a subset of $\der_k(A)$ and $\tau: \g \to \der_k(A)$ is
the inclusion map, the pair $(\g, \tau)$ is a Lie-Rinehart algebra if
and only if $\g$ is closed under the $A$-module and $k$-Lie algebra
structures of $\der_k(A)$. We are mainly interested in Lie-Rinehart
algebras of this type, and usually omit $\tau$ from the notation.

Let $\g$ be a Lie-Rinehart algebra. For any $A$-module $M$, we define a
\emph{$\g$-connection} on $M$ to be an $A$-linear map $\nabla: \g \to
\enm_k(M)$ such that
\begin{equation} \label{e:dp}
\nabla_D(am) = a \nabla_D(m) + D(a) \; m
\end{equation}
for all $D \in \g, \; a \in A, \; m \in M$. We say that $\nabla$
satisfies the \emph{derivation property} when condition (\ref{e:dp})
holds for all $D \in \g$. If $\nabla: \g \to \enm_k(M)$ is a $k$-linear
map that satisfies the derivation property, we call $\nabla$ a
$k$-linear $\g$-connection on $M$. A \emph{connection} on $M$ is a
$\g$-connection on $M$ with $\g = \der_k(A)$.

Let $\nabla$ be a $\g$-connection on $M$. We define the
\emph{curvature} of $\nabla$ to be the $A$-linear map $R_{\nabla}: \g
\wedge \g \to \enm_A(M)$ given by
\begin{equation*}
R_{\nabla}(D \wedge D') = [ \nabla_D, \nabla_{D'} ] - \nabla_{[D,D']}
\end{equation*}
for all $D,D' \in \g$. We say that $\nabla$ is an \emph{integrable
$\g$-connection} if $R_{\nabla} = 0$.

We define $\mc(A,\g)$ to be the category of \emph{modules with
$\g$-connections}. The objects in $\mc(A,\g)$ are pairs $(M,\nabla)$,
where $M$ is an $A$-module and $\nabla$ is a $\g$-connection on $M$,
and the morphisms $\phi: (M,\nabla) \to (M', \nabla')$ in $\mc(A,\g)$
are the \emph{horizontal maps}, i.e. $A$-linear homomorphisms $\phi: M
\to M'$ such that $\phi \nabla_D = \nabla'_D \phi$ for all $D \in \g$.
The category $\mc(A,\g)$ is an Abelian $k$-category, and we write
$\mc(A) = \mc(A,\g)$ when $\g = \der_k(A)$.

Let $\mic(A,\g)$ be the full subcategory of $\mc(A,\g)$ of
\emph{modules with integrable $\g$-connections}. This is an Abelian
subcategory with many nice properties. In fact, there is an associative
$k$-algebra $\Delta(A,\g)$ such that the category $\mic(A,\g)$ is
equivalent to the category of left modules over $\Delta(A,\g)$. When
$\g \subseteq \der_k(A)$, $\Delta(A,\g)$ is the subalgebra of
$\diff(A)$ generated by $A$ and $\g$, where $\diff(A)$ denotes the ring
of differential operators on $A$ in the sense of Grothendieck
\cite{gr67}. When $\g = \der_k(A)$, the algebra $\Delta(A) =
\Delta(A,\g)$ is called the \emph{derivation algebra}.

We recall that when $A$ is a regular $k$-algebra, a \emph{connection}
on $M$ is usually defined as a $k$-linear map $\nabla: M \to M
\otimes_A \kdf A$ such that $\nabla(am) = a \nabla(m) + m \otimes d(a)$
for all $a \in A, \; m \in M$, see Katz \cite{ka70}. Moreover, the
\emph{curvature} of $\nabla$ is usually defined as the $A$-linear map
$R_{\nabla}: M \to M \otimes_A \kdf{A}^2$ given by $R_{\nabla} =
\nabla^1 \circ \nabla$, where $\nabla^1$ is the natural extension of
$\nabla$ to $M \otimes \kdf A$, and $\nabla$ is an \emph{integrable
connection} if $R_\nabla = 0$.

Let $A$ be any commutative $k$-algebra. For expository purposes, we
define an \emph{$\Omega$-connection} on an $A$-module $M$ to be a
connection on $M$ in the sense of the preceding paragraph. We define
$\omc(A)$ to be the Abelian $k$-category of modules with
$\Omega$-connections, and $\omic(A)$ to be the full Abelian subcategory
of modules with integrable $\Omega$-connections.

\begin{lem} \label{l:mcp}
Let $A$ be a regular local $k$-algebra essentially of finite type, and
let $M$ be a finitely generated $A$-module. If there is an
$\Omega$-connection on $M$, then $M$ is free.
\end{lem}
\begin{proof}
Let $\{ t_1, \dots, t_d \}$ be a set of regular parameters of $A$, and
let $\delta_i$ be the derivation on $k[t_1, \dots, t_d]$ such that
$\delta_i(t_j) = \delta_{ij}$ for $1 \le i \le d$. If $A$ is
essentially of finite type over $k$, then $\delta_i$ extends to a
derivation of $A$. This implies that any module $M$ that admits an
$\Omega$-connection is free, see for instance Borel et al. \cite{bo87},
section VI, proposition 1.7.
\end{proof}

\begin{lem} \label{l:eq-mc}
There is a natural functor $\omc(A) \to \mc(A)$, and an induced functor
$\omic(A) \to \mic(A)$. If $\kdf A$ and $\der_k(A)$ are projective
$A$-modules of finite presentation, then these functors are
equivalences of categories.
\end{lem}
\begin{proof}
Any $\Omega$-connection on $M$ induces a connection on $M$, and this
assignment preserves integrability. Moreover, any connection $\nabla$
on $M$ may be considered as a $k$-linear map $M \to
\hmm_A(\der_k(A),M)$, given by $m \mapsto \{ D \mapsto \nabla_D(m) \}$.
It is sufficient to show that the natural map $M \otimes_A \kdf A \to
\hmm_A(\der_k(A),M)$, given by $m \otimes \omega \mapsto \{ D \mapsto
\phi_D(\omega) m \}$, is an isomorphism. But this is clearly the case
when $\kdf A$ and $\der_k(A)$ are projective $A$-modules of finite
presentation.
\end{proof}

We see that if $A$ is a regular $k$-algebra essentially of finite type,
then there is a bijective correspondence between (integrable)
connections on $M$ and (integrable) $\Omega$-connections on $M$ for any
$A$-module $M$. In contrast, there are many modules that admit
connections but not $\Omega$-connections when $A$ is a singular
$k$-algebra.

When $A$ is a complete local $k$-algebra, we must replace the tensor
and wedge products with their formal analogues. We remark that with
this modification, lemma \ref{l:mcp} and lemma \ref{l:eq-mc} hold
for any complete local Noetherian $k$-algebra $A$ with residue field
$k$.

\section{Elementary properties of connections}

Let $k$ be an algebraically closed field of characteristic $0$, let
$A$ be a commutative $k$-algebra, and let $\g$ be a Lie-Rinehart
algebra. If $(M,\nabla)$ and $(M', \nabla')$ are modules with
$\g$-connections, then there are natural induced $\g$-connections on
the $A$-modules $M \oplus M'$ and $\hmm_A(M,M')$. These
$\g$-connections are integrable if and only if $\nabla$ and
$\nabla'$ are integrable $\g$-connections.

\begin{lem} \label{l:dirsum}
For any $A$-modules $M,M'$, $M \oplus M'$ admits a $\g$-connection if
and only if $M$ and $M'$ admit $\g$-connections.
\end{lem}

\begin{lem} \label{l:dual}
For any reflexive $A$-module $M$, $M^\vee = \hmm_A(M,A)$ admits a
$\g$-connection if and only if $M$ admits a $\g$-connection.
\end{lem}

In view of lemma \ref{l:mcp} and lemma \ref{l:eq-mc}, we expect that
any $A$-module that admits a connection must be locally free outside
the singular locus of $\spec(A)$. In order to prove this, we need some
results on localizations of connections.

Let $A \to \loc A$ be the localization given by a multiplicatively
closed subset $S \subseteq A$. Since any derivation of $A$ can be
extended to a derivation of $\loc A$, we see that $\loc \g$ is a
Lie-Rinehart algebra of $\loc A/k$.

\begin{lem} \label{l:loc-conn}
Localization gives a functor $\mc(A,\g) \to \mc(\loc A,\loc \g)$ and an
induced functor $\mic(A,\g) \to \mic(\loc A, \loc \g)$ for any
multiplicatively closed subset $S \subseteq A$.
\end{lem}

Let $A \to \widehat A$ be the $m$-adic completion of $A$ given by a
maximal ideal $m \subseteq A$. Since any derivation of $A$ can be
extended to a derivation of $\widehat A$, we see that $\widehat A
\otimes_A \g$ is a Lie-Rinehart algebra of $\widehat A/k$. Moreover, if
$A$ is Noetherian and $\g$ is a finitely generated $A$-module, then
$\widehat A \otimes_A \g \cong \widehat \g$.

\begin{lem} \label{l:compl-conn}
If $A$ is a Noetherian $k$-algebra and $\g$ is a finitely generated
$A$-module, $m$-adic completion gives a functor $\mc(A,\g) \to
\mc(\widehat A,\widehat \g)$ and an induced functor $\mic(A,\g) \to
\mic(\widehat A, \widehat \g)$ for any maximal ideal $m \subseteq A$.
\end{lem}

In particular, if $A$ is essentially of finite type over $k$, then
there are localization functors $\mc(A) \to \mc(\loc A)$ and $\mic(A)
\to \mic(\loc A)$ for any multiplicatively closed subset $S \subseteq
A$, and $m$-adic completion functors $\mc(A) \to \mc(\widehat A)$ and
$\mic(A) \to \mic(\widehat A)$ for any maximal ideal $m \subseteq A$.

\begin{lem} \label{l:conn-reg}
Let $A$ be a $k$-algebra essentially of finite type, and let $M$ be a
finitely generated $A$-module. If there is a connection on $M$, then
$M_p$ is a locally free $A_p$-module for all prime ideals $p \subseteq
A$ such that $A_p$ is a regular local ring.
\end{lem}

This lemma also holds when $A$ is a complete local Noetherian
$k$-algebra with residue field $k$. We remark that lemma
\ref{l:conn-reg} gives a necessary condition for a module to admit
connections, and it is well-known that maximal Cohen-Macaulay modules
satisfy this condition.

\begin{lem} \label{l:conn-is}
Let $A$ be a $k$-algebra essentially of finite type, let $\g$ be a
Lie-Rinehart algebra of $A/k$ that is finitely generated as an
$A$-module, and let $M$ be a finitely generated $A$-module. For any
maximal ideal $m \subseteq A$, we write $\widehat \g$ and $\widehat M$
for the $m$-adic completions of $\g$ and $M$, and consider the
following statements:
\begin{enumerate}
\item $M$ admits a $\g$-connection
\item $M_m$ admits a $\g_m$-connection
\item $\widehat M$ admits a $\widehat \g$-connection
\end{enumerate}
Then we have $(1) \Rightarrow (2) \Leftrightarrow (3)$. Moreover, if
$M_p$ admits a $\g_p$-connection for all prime ideals $p \neq m$ in
$A$, then $(1) \Leftrightarrow (2)$.
\end{lem}
\begin{proof}
The implications $(1) \Rightarrow (2) \Rightarrow (3)$ is a direct
consequence of lemma \ref{l:loc-conn} and lemma \ref{l:compl-conn}. By
the obstruction theory for connections, see Eriksen, Gustavsen
\cite{er-gu06}, it follows that $(2) \Leftrightarrow (3)$. Furthermore,
if there is a $\g_p$-connection on $M_p$ for all prime ideals $p \neq
m$, it also follows that $(3)$ implies $(1)$.
\end{proof}

\section{Graded connections}

Let $k$ be an algebraically closed field of characteristic $0$, and let
$A$ be a \emph{quasi-homogeneous} $k$-algebra, i.e. a positively graded
$k$-algebra of the form $A \cong S/I$, where $S = k[x_1, \dots, x_n]$
is a graded polynomial ring with $\deg(x_i) > 0$ for $1 \le i \le n$,
and $I$ is a homogeneous ideal in $S$.

We see that $\der_k(A)$ has a natural grading induced by the grading of
$A$ such that the homogeneous derivations $D \in \der_k(A)$ of degree
$\omega$ satisfy $D(A_i) \subseteq A_{i+\omega}$ for all integers $i
\ge 0$. We also notice that $\der_k(A)$ is a graded $k$-Lie algebra,
i.e. $[\der_k(A)_i , \der_k(A)_j] \subseteq \der_k(A)_{i+j}$ for all
integers $i,j$.

We say that a Lie-Rinehart algebra $(\g, \tau)$ is \emph{graded} if
$\g$ is a graded $A$-module and $k$-Lie algebra, and $\tau: \g \to
\der_k(A)$ is a graded homomorphism of $A$-modules and $k$-Lie
algebras.

Let $\g$ be a graded Lie-Rinehart algebra. For any graded $A$-module
$M$, we define a \emph{graded $\g$-connection} on $M$ to be a
$\g$-connection $\nabla$ on $M$ such that $\nabla_D(M_i) \subseteq
M_{i+\omega}$ for any integer $i$ and for any homogeneous element $D
\in \g$ of degree $\omega$.

\begin{lem} \label{l:gr-conn}
Let $A$ be a quasi-homogeneous $k$-algebra, let $\g$ be a graded
Lie-Rinehart algebra of $A/k$ that is finitely generated as an
$A$-module, and let $m$ be the graded maximal ideal of $A$. For any
finitely generated graded $A$-module $M$, the following conditions are
equivalent:
\begin{enumerate}
\item $M$ admits a graded $\g$-connection
\item $M_m$ admits a $\g_m$-connection
\end{enumerate}
\end{lem}
\begin{proof}
The functor $M \mapsto M_m$ is faithfully exact on the category of
finitely generated graded $A$-modules since $A$ has a unique graded
maximal ideal $m$. Hence the result follows as in lemma
\ref{l:conn-is}.
\end{proof}

Let $A$ be a quasi-homogeneous $k$-algebra, and let $M$ be a finitely
generated graded $A$-module. We consider a graded presentation of $M$
of the form
\begin{equation} \label{e:modpres}
    0 \gets M \xleftarrow{\rho} L_0 \xleftarrow{d_0} L_1,
\end{equation}
where $L_0$ and $L_1$ are free graded $A$-modules of finite rank, with
homogeneous $A$-linear bases $\{ e_i \}$ and $\{ f_i \}$ respectively,
and $d_0$ is a graded $A$-linear homomorphism. If we write $( a_{ij} )$
for the matrix of $d_0$ with respect to the chosen bases, then $a_{ij}$
is homogeneous with $\deg(a_{ij}) = \deg(e_i) - \deg(f_j)$ for $1 \le i
\le \rk(L_0), \; 1 \le j \le \rk(L_1)$.

Let $D \in \der_k(A)$ be a homogeneous derivation of degree $\omega$.
Then there is a natural action of $D$ on $L_0$ and $L_1$, i.e. $D(a
e_i) = D(a) e_i$ and $D(a f_j) = D(a) f_j$ for any $a \in A$ and any
$e_i \in L_0, \; f_j \in L_1$. For simplicity, we shall denote the
induced $k$-linear endomorphisms by $D: L_n \to L_n$ for $n = 0,1$, and
write $D(d_0): L_1 \to L_0$ for the $A$-linear homomorphism given by
$D(d_0) = D d_0 - d_0 D$. Notice that $D: L_n \to L_n$ is graded of
degree $\omega$, hence $D(d_0)$ is also graded of degree $\omega$.

\begin{lem} \label{l:connpres}
Let $D \in \der_k(A)_w$ for some integer $\omega$, and let $\nabla_D
\in \enm_k(M)$ be a $k$-linear endomorphism with derivation property
with respect to $D$ such that $\nabla_D(M_i) \subseteq M_{i+\omega}$
for all integers $i$. Then there exist $A$-linear endomorphisms $\phi_D
\in \enm_A(L_0)_\omega$ and $\psi_D \in \enm_A(L_1)_\omega$ with
$D(d_0) = d_0 \psi_D - \phi_D d_0$ such that $\nabla_D$ is induced by
$D + \phi_D: L_0 \to L_0$.
\end{lem}
\begin{proof}
Consider the map $\nabla_D \, \rho - \rho D: L_0 \to M$, and notice
that it is a graded $A$-linear homomorphism of degree $\omega$. Hence
we can find a graded $A$-linear endomorphism $\phi_D: L_0 \to L_0$ of
degree $\omega$ such that $\nabla_D \, \rho = \rho (D + \phi_D)$, and
this implies that $(D+\phi_D) d_0(x) \in \im(d_0)$ for all $x \in L_1$.
So we can find a graded $A$-linear homomorphism $\psi_D: L_1 \to L_1$
of degree $\omega$ such that $D(d_0) = d_0 \psi_D - \phi_D d_0$.
\end{proof}

If we write $( p_{ij} )$ for the matrix of $\phi_D$ and $( c_{ij} )$
for the matrix of $\psi_D$ with respect to the chosen bases, then
$p_{ij}$ is homogeneous with $\deg(p_{ij}) = \deg(e_j) - \deg(e_i) +
\omega$ for $1 \le i,j \le \rk(L_0)$ and $c_{ij}$ is homogeneous with
$\deg(c_{ij}) = \deg(f_j) - \deg(f_i) + \omega$ for $1 \le i,j \le
\rk(L_1)$.

\section{Maximal Cohen-Macaulay modules}

Let $k$ be an algebraically closed field of characteristic $0$, and
let $A$ be a complete local Noetherian $k$-algebra with residue
field $k$. We say that a finitely generated $A$-module $M$ is
\emph{maximal Cohen-Macaulay} (MCM) if $\depth(M) = \dim(A)$, and
that $A$ is a \emph{Cohen-Macaulay} (CM) ring if $A$ is MCM as an
$A$-module. In the rest of this paper, we shall assume that $A$ is a
CM ring and an isolated singularity.

We say that $A$ has \emph{finite CM representation type} if the number
of isomorphism classes of indecomposable MCM $A$-modules is finite. The
singularities of finite CM representation type has been classified when
$A$ has dimension $d \le 2$, but not in higher dimensions.

Let $A$ be the complete local ring of a \emph{simple hypersurface
singularity}, see Arnold \cite{arn81} and Wall \cite{wal84}. Then $A
\cong k[[z_0, \dots, z_d]] / (f)$, where $d \ge 1$ is the dimension
of $A$ and $f$ is of the form
\begin{align*}
A_n: \; & f = z_0^2 + z_1^{n+1} + z_2^2 + \dots + z_d^2 & n \ge 1\\
D_n: \; & f = z_0^2 z_1 + z_1^{n-1} + z_2^2 + \dots + z_d^2 & n \ge 4\\
E_6: \; & f = z_0^3 + z_1^4 + z_2^2 + \dots + z_d^2 & \\
E_7: \; & f = z_0^3 + z_0 z_1^3 + z_2^2 + \dots + z_d^2 & \\
E_8: \; & f = z_0^3 + z_1^5 + z_2^2 + \dots + z_d^2 &
\end{align*}
The simple singularities are exactly the hypersurface singularities of
finite CM repre\-sentation type, see Kn\"orrer \cite{kno87} and
Buchweitz, Greuel, Schreyer \cite{bu-gr-sc87}. Moreover, if $A$ is
Gorenstein and of finite CM representation type, then $A$ is a simple
singularity, see Herzog \cite{her78}.

Assume that $A$ is a hypersurface $A = S/(f)$, where $S$ is a power
series $k$-algebra. A \emph{matrix factorization} of $f$ is a pair
$(\phi, \psi)$ of square matrices with entries in $S$ such that $\phi
\psi = \psi \phi = I$. We say that $(\phi, \psi)$ is a \emph{reduced
matrix factorization} if the entries in $\phi$ and $\psi$ are non-units
in $S$. By Eisenbud \cite{eis80}, there is a bijective correspondence
between reduced matrix factorizations of $f$ and MCM $A$-modules
without free summands, given by the assignment $(\phi,\psi) \mapsto
\coker(\phi)$.

Let $f \in S = k[[z_0, \dots, z_d]]$ be the equation of the simple
singularity $A = S/(f)$ of dimension $d$, and let $(\phi, \psi)$ be a
reduced matrix factorization of $f \in S$. Then the pair $(\phi',
\psi')$ given by
    \[ \phi' = \begin{pmatrix} u I & - \psi \\ \phi & v I
    \end{pmatrix}, \psi' = \begin{pmatrix} v I & \psi \\ - \phi & u I
    \end{pmatrix} \]
is a reduced matrix factorization of $f' = f + uv \in S' = S[[u,v]]$,
and $A' = S'/(f')$ is isomorphic to the simple singularity of dimension
$d+2$ corresponding to $A$. By Kn\"orrer's periodicity theorem, this
assignment induces a bijective correspondence between MCM $A$-modules
without free summands and MCM $A'$-modules without free summands, see
Kn\"orrer \cite{kno87} and Schreyer \cite{sch87}.

A complete list of indecomposable MCM modules over simple curve
singularities was given in Greuel, Kn\"orrer \cite{gr-kn85}, and the
corresponding matrix factorizations were given in Eriksen \cite{er06}
and Yoshino \cite{yos90}. A complete list of indecomposable MCM modules
over simple surface singularities can be obtained from the irreducible
representations of the finite subgroups of $\slg(2,k)$. In higher
dimensions, a complete list of indecomposable MCM modules over simple
singularities can be obtained using Kn\"orrer's periodicity theorem.

\begin{lem} \label{l:mfid}
Let $(\phi, \psi)$ be a reduced matrix factorization of $f \in S$, let
$A = S/(f)$, and let $M = \coker(\phi)$. Then we have:
\begin{enumerate}
\item $\coker(\phi^t) \cong M^\vee$, the $A$-linear dual of $M$,
\item $\coker(\psi) \cong \syz^1(M)$, the first reduced syzygy of $M$.
\end{enumerate}
\end{lem}

We remark that by lemma \ref{l:dual}, there is a connection on $M$ if
and only if there is a connection on $M^\vee$. It is not difficult to
prove that there is a $k$-linear connection on $M$ if and only if there
is a $k$-linear connection on $\syz^1(M)$ using the above result, and a
different proof of this fact was given in K\"allstr\"om \cite{kal05},
theorem 2.2.8. However, we do not know if is true that there is a
connection on $M$ if and only if there is a connection of $\syz^1(M)$.

\section{Connections on MCM modules}
\label{s:mcm-conn}

Let $k$ be an algebraically closed field of characteristic $0$, and let
$A$ be a complete local CM $k$-algebra with residue field $k$ that is
an isolated singularity. In this section, we study the existence of
connections on MCM $A$-modules under these conditions. We focus on the
cases when $A$ has finite CM representation type.

\subsection{Dimension zero}

When $A$ is a complete local ring of a zero-dimensional singularity, it
has finite CM representation type if and only if it is of the form $A =
k[[x]]/(x^{n+1})$ for some integer $n \ge 1$, see Herzog \cite{her78},
satz 1.5. We consider these singularities as the zero-dimensional
simple singularities of type $A_n$.

We remark that there are $n-1$ reduced matrix factorizations of $x^n$,
given by $x^n = x^i \cdot x^{n-i}$ for $1 \le i \le n-1$. Since the
natural action of the Euler derivation $E = x \frac{\partial}{\partial
x}$ on $x^i$ is given by $E(x^i) = i x^i$, we see that any MCM
$A$-module admits an integrable connection.

\subsection{Dimension one} \label{ss:dim1}

When $A$ is the complete local ring of a simple curve singularity, it
was shown in Eriksen \cite{er06} that any MCM $A$-module admits an
integrable connection.

Let $A$ be the complete local ring of a curve singularity. We say that
a local ring $B$ \emph{birationally dominates} $A$ if $A \subseteq B
\subseteq A^*$, where $A^*$ is the integral closure of $A$ in its total
quotient ring. It is known that $A$ has finite CM representation type
if and only if it birationally dominates the complete local ring of a
simple curve singularity, see Greuel, Kn\"orrer \cite{gr-kn85}.

This result leads to a complete classification of curve
singularities of finite CM representation type. The non-Gorenstein
curve singularities of finite CM representation type are of the
following form:
\begin{align*}
D^s_n: \; & A = k[[x,y,z]]/(x^2-y^n,xz,yz) \text{ for } n \ge 2 \\
E^s_6: \; & A = k[[t^3,t^4,t^5]] \subseteq k[[t]] \\
E^s_7: \; & A = k[[x,y,z]]/(x^3-y^4,xz-y^2,y^2z-x^2,yz^2-xy) \\
E^s_8: \; & A = k[[t^3,t^5,t^7]] \subseteq k[[t]]
\end{align*}
Using \textsc{Singular} \cite{gps05} and our library {\sc conn.lib}
\cite{er-gu05}, we show that not all MCM $A$-modules admit connections
in these cases. In fact, the canonical module $\omega_A$ does not admit
a connection when $A$ is the complete local ring of the singularities
$E^s_6, E^s_7, E^s_8$ or $D^s_n$ for $n \le 100$.

The \emph{monomial curve singularities} is an interesting class of
curve singularities not necessarily of finite CM representation type.
Let $\Gamma \subseteq \n$ be a \emph{numerical semigroup}, i.e. a
sub-semigroup $\Gamma \subseteq \n$ such that the complement $H = \n
\setminus \Gamma$ is finite, and let $A = k[[\Gamma]]$ be the complete
local ring of the corresponding monomial curve singularity, i.e. the
subalgebra $k[[t^{a_1}, t^{a_2}, \dots, t^{a_r}]] \subseteq k[[t]]$,
where $\{ a_1, a_2, \dots, a_r \}$ is the minimal set of generators of
$\Gamma$. Let the \emph{Frobenius number} $g$ of $\Gamma$ be the
maximal element of $H$. It is well-known that $A$ is Gorenstein if and
only if $\Gamma$ is \emph{symmetric}, i.e. $a \in \Gm$ if and only if
$g-a \not \in \Gamma$ for any integer $a \in \z$. If $\Gamma$ is not
symmetric, we denote by $\Delta$ the non-empty set $\Delta = \{ h \in
H: g-h \in H \}$.

Let $\Lambda$ be a set such that $\Gamma \subseteq \Lambda \subseteq
\n$ and $\Gamma + \Lambda \subseteq \Lambda$, and consider the module
$M = k[[\Lambda]]$ with $k$-linear basis $\{ t^a: a \in \Lambda \}$ and
the obvious action of $A$. It is clear that $M = k[[\Lambda]]$ is an
MCM $A$-module of rank one, and that $M = k[[\Lambda]]$ is isomorphic
to $M' = k[[\Lambda']]$ if and only if $\Lambda = \Lambda'$. In fact,
one may show that the rank one MCM modules over $A=k[[\Gamma]]$ of the
form $M = k[[\Lambda]]$ are exactly the \emph{gradable} ones, i.e. the
$A$-modules induced by graded modules over the quasi-homogeneous
algebra $k[\Gamma] = k[t^{a_1},t^{a_2}, \dots, t^{a_r}] \subseteq
k[t]$.

\begin{prop} \label{p:nconn-mon} Let $A = k[[\Gamma]]$ be a monomial
curve singularity, let $\Lambda$ be set such that $\Gamma \subseteq
\Lambda \subseteq \n$ and $\Gamma + \Lambda \subseteq \Lambda$, and let
$M = k[[\Lambda]]$ be the corresponding rank one MCM $A$-module. If
$g-s \in \Lambda$ and $g \not \in \Lambda$, where $s = \max \Delta$,
then $M$ does not admit connections.
\end{prop}
\begin{proof}
We see that $E = t \frac{\partial}{\partial t}$, $D =t^g E$ and $D' =
t^s E$ are derivations of $A$, see Eriksen \cite{er03}, the remarks
preceding lemma 8. If $\nabla$ is a connection on $M$, then there
exists an element $f = f_0 + f_+ \in k[[\Lambda]]$ with $f_0 \in k$ and
$f_+ \in (t)$ such that $\nabla_E(t^\lambda) = E(t^\lambda) + f
t^\lambda$ for all $\lambda \in \Lambda$. Since $M$ is torsion free, we
must have $\nabla_D(t^\lambda) = D(t^\lambda) + t^g f t^\lambda \in M$
for all $\lambda \in \Lambda$. For $\lambda = 0$, this condition
implies that $f_0 = 0$. Similarly, we must have $\nabla_{D'}(t^\lambda)
= D'(t^\lambda) + t^s f t^\lambda \in M$ for all $\lambda \in \Lambda$.
For $\lambda = g-s$, this condition implies that $(g-s) + f_0 = 0$,
which is a contradiction.
\end{proof}

\begin{thm} \label{t:conn-mon}
Let $A$ be the complete local ring of a monomial curve singularity.
Then all gradable MCM $A$-module of rank one admits a connection if and
only if $A$ is Gorenstein.
\end{thm}
\begin{proof}
Let $A = k[[\Gamma]]$  for a numerical semigroup $\Gamma$, and $M =
k[[\Lambda]]$ for a set $\Lambda$ with $\Gamma \subseteq \Lambda
\subseteq \n$ and $\Gamma + \Lambda \subseteq \Lambda$. If $A$ is
Gorenstein, then $\der_k(A)$ is generated by the Euler derivation $E =
t \frac{\partial}{\partial t}$ and the trivial derivation $D = t^g E$,
see Eriksen \cite{er03}, lemma 8 and the following remarks. Hence the
natural action of $E$ and $D$ on $M$ induces a connection on $M$ since
$g + (\Lambda \setminus \{ 0 \}) \subseteq \Gamma$ and $D(1) = 0$. On
the other hand, if $A$ is not Gorenstein, then the set $\Delta$ is
non-empty, and $\Lambda = \Gamma \cup (g - s + \Gamma)$ satisfies the
conditions $\Gamma \subseteq \Lambda \subseteq \n$ and $\Gamma +
\Lambda \subseteq \Lambda$. Since $g \not \in \Lambda$, it follows from
proposition \ref{p:nconn-mon} that $M$ does not admit connections.
\end{proof}

Let us consider the non-Gorenstein monomial curve singularities $E^s_6$
and $E^s_8$ of finite CM representation type. By Yoshino \cite{yos90},
theorem 15.14, all MCM $A$-modules are gradable in these cases. For
$E^s_6$, we have $H = \{ 1,2 \}$, and the possibilities for $\Lambda$
(with $\Lambda \neq \Gamma$) are
    \[ \Lambda_1 = \Gamma \cup \{ 1 \}, \quad \Lambda_2 = \Gamma \cup \{
    2 \}, \quad \Lambda_{12} = \Gamma \cup \{ 1, 2 \} \]
The corresponding modules $M_1, M_2, M_{12}$ are the non-free rank one
MCM $A$-modules. Using the method from the proof of proposition
\ref{p:nconn-mon}, we see that the modules $M_2$ and $M_{12}$ admit
connections, while $M_1$ does not. One may show that $M_1$ is the
canonical module in this case. A similar consideration for $E^s_8$
shows that $M_{14}$, $M_2$, $M_4$, $M_{24}$ and $M_{124}$ are the
non-free rank one MCM $A$-modules, and that the canonical module $M_2$
is the only rank one MCM $A$-module that does not admit connections.

Let us also consider a non-Gorenstein monomial curve singularity not of
finite CM representation type. When $A = k[[t^4, t^5, t^6, t^7]]$, we
have $H = \{ 1,2,3 \}$, and the possibilities for $\Lambda$ (with
$\Lambda \neq \Gamma$) are
\begin{multline*}
\Lambda_{1} = \Gamma \cup \{ 1 \}, \quad \Lambda_2 = \Gamma \cup \{ 2
\}, \quad \Lambda_3 = \Gamma \cup \{ 3 \}, \quad \Lambda_{12} = \Gamma
\cup \{ 1,2 \}, \\
\Lambda_{13} = \Gamma \cup \{ 1,3 \}, \quad \Lambda_{23} = \Gamma \cup
\{ 2,3 \}, \quad \Lambda_{123} = \Gamma \cup \{ 1, 2, 3 \}
\end{multline*}
The corresponding modules $M_1, M_2, M_3$, $M_{12}, M_{13}, M_{23}$
and $M_{123}$ are the non-free rank one gradable MCM $A$-modules.
One may show that the modules $M_3$, $M_{23}$ and $M_{123}$ admit
connections, that the canonical module $M_{12}$ and the module
$M_{13}$ admit $k$-linear connections but not connections, while the
modules $M_1$ and $M_2$ do not even admit $k$-linear connections.

Finally, we remark that any connection on an MCM $A$-module is
integrable when $A$ is a monomial curve singularity.

\subsection{Dimension two}

When $A$ is the complete local ring of a surface singularity, it has
finite CM representation type if and only if it is a quotient
singularity of the form $A = S^G$, where $S = k[[x,y]]$ and $G$ is a
finite subgroup of $\glg(2,k)$ without pseudo-reflections. This fact
was proven independently in Auslander \cite{aus86} and Esnault
\cite{esn85}. Moreover, there is a bijective correspondence between MCM
$A$-modules and finite dimensional representations of $G$.

It is not difficult to see that any MCM module over a quotient surface
singularity admits an integrable connection, and Jan Christophersen was
the first to point this out to us.

\begin{prop} \label{p:conn-quot}
Let $A$ be the complete local ring of a quotient singularity $A = S^G$,
where $S = k[[x_1, \dots, x_n]]$ and $G \subseteq \glg(n,k)$ is a
finite subgroup without pseudo-reflections. For any finite dimensional
representation $\rho: G \to \enm_k(V)$, the MCM $A$-module $M = (S
\otimes_k V)^G$ admits an integrable connection.
\end{prop}
\begin{proof}
There is a canonical integrable connection $\nabla': \der_k(S) \to
\enm_k(S \otimes_k V)$ on the free $S$-module $S \otimes_k V$, given by
    \[ \nabla'_D(\sum s_i \otimes v_i) = \sum D(s_i) \otimes v_i \]
for any $D \in \der_k(S), \; s_i \in S, \; v_i \in V$. But the natural
map $\der_k(S)^G \to \der_k(S^G)$ is an isomorphism, see Kantor
\cite{kan71} or Schlessinger \cite{sch71}, hence $\nabla'$ induces an
integrable connection $\nabla$ on $M$.
\end{proof}

\begin{thm} \label{t:quotsing-dim2}
Let $A$ be the complete local ring of a surface singularity of finite
CM representation type. Then any MCM $A$-module admits an integrable
connection.
\end{thm}
\begin{proof}
By the comments preceding proposition \ref{p:conn-quot}, we may assume
that $A = S^G$, where $S = k[[x,y]]$ and $G$ is a finite subgroup of
$\glg(2,k)$ without pseudo-reflections, and that $M = (S \otimes_k
V)^G$ for a finite dimensional representation $\rho: G \to \enm_k(V)$
of $G$. Hence $M$ admits an integrable connection by proposition
\ref{p:conn-quot}.
\end{proof}

We recall that the simple surface singularities are precisely the
quotient surface singularities that are Gorenstein. We may also
characterize them as the quotient surface singularities $A = S^G$ where
$G$ is a subgroup of $\slg(2,k)$. In particular, we see that any MCM
$A$-module over a simple surface singularity admits an integrable
connection.

Let $A$ be the complete local ring of any surface singularity. We
remark that theorem \ref{t:quotsing-dim2} can be generalized as
follows: Let us consider a finite Galois extension $L$ of $K$, where
$K$ is the field of fractions of $A$, and let $B$ be the integral
closure of $A$ in $L$. If the extension $A \subseteq B$ is unramified
at all height one prime ideals, we say that it is a \emph{Galois
extension}. In Gustavsen, Ile \cite{gu-il06a}, it was proven that if
$M$ is an MCM $A$-module such that $(M \otimes_A B)^{\vee\vee}$ is a
free $B$-module for some Galois extension $A \subseteq B$, then $M$
admits an integrable connection. In particular, if $A$ is a rational
surface singularity, it follows that any rank one MCM $A$-module admits
an integrable connection.

If $A = S^G$ is a quotient singularity, then $A \subseteq S = k[[x,y]]$
is a Galois extension such that $(M \otimes_A S)^{\vee\vee}$ is free
for any MCM $A$-module $M$. In contrast, if $A$ is not a quotient
singularity, there exists an MCM $A$-module $M$ such that $(M \otimes_A
B)^{\vee\vee}$ is non-free for any Galois extension $A \subseteq B$,
see Gustavsen, Ile \cite{gu-il06a}. Nevertheless, $M$ may still admit
an integrable connection. In fact, any MCM module over a simple
elliptic surface singularity admits an integrable connection, see Kahn
\cite{ka88}, and this result has been generalized to quotients of
simple elliptic surface singularities in Gustavsen, Ile
\cite{gu-il06b}. Kurt Behnke has pointed out that it might be true for
cusp singularities as well, see Behnke \cite{be89}. More generally, it
is probable that any MCM module over a log canonical surface
singularity admits an integrable connection.

When $A$ is a surface singularity, we have not found any examples of an
MCM module that does not admit an integrable connection.

\subsection{Higher dimensions}

The main result of this paper is that when $A$ is the complete local
ring of a simple singularity of dimension $d \ge 3$, an MCM $A$-module
$M$ admits a connection only if $M$ is free. In Eriksen, Gustavsen
\cite{er-gu06}, we used \textsc{Singular} \cite{gps05} and our library
{\sc conn.lib} \cite{er-gu05} to prove this result when $A$ is a simple
threefold singularity of type $A_n$, $D_n$ or $E_n$ with $n \le 50$.
Using graded techniques, we are able to prove this result for any
simple singularity of dimension $d \ge 3$:

\begin{thm} \label{t:simple-dim3}
Let $A$ be the complete local ring of a simple singularity of dimension
$d \ge 3$. Then an MCM $A$-module $M$ admits a connection if and only
if $M$ is free.
\end{thm}
\begin{proof}
We claim that when $A$ is a simple threefold singularity, $\g$ is the
submodule of $\der_k(A)$ generated by the trivial derivations, and $M$
is a non-free MCM $A$-module, then $M$ does not admit $\g$-connections.

Let us first prove that the claim implies the theorem. If $A$ is a
simple threefold singularity and $M$ is a non-free MCM $A$-module that
admits a connection, then its restriction to $\g$ is a $\g$-connection
on $M$. So we may assume that $A$ has dimension $d > 3$. Assume that
there is a non-free MCM $A$-module $M$ that admits a connection.
Clearly, $A = S[[u_1, \dots, u_{d-3}]]/(f)$ with $f = f_0 + u_1^2 +
\dots + u_{d-3}^2$, where $S/(f_0)$ is the simple threefold singularity
of the same type as $A$. Moreover, $M$ is given by a reduced matrix
factorization $(\phi, \psi)$ of $f$. Let us write $\g$ for the
submodule of $\der_k(S/(f_0))$ generated by the trivial derivations.
Since $A/(u_1, \dots, u_{d-3}) \cong S/(f_0)$, $M_0 = M/(u_1, \dots,
u_{d-3})M$ is a non-free MCM $S/(f_0)$-module, given by the reduced
matrix factorization $(\phi \otimes_A A/(u_1, \dots, u_{d-3}), \psi
\otimes_A A/(u_1, \dots, u_{d-3}) )$ of $f_0$. Since $M$ admits a
connection, its restriction to $\g_A = A \cdot \g \subseteq \der_k(A)$,
the $A$-submodule of $\der_k(A)$ generated by the trivial derivations
of $S/(f)$, is a $\g_A$-connection on $M$. Since any $\g_A$-connection
on $M$ induces a $\g$-connection on $M_0$, this contradicts the claim.

It remains to prove the claim. Let $A' = k[[z_0,z_1,z_2,z_3]]/(f)$ be
the complete local ring of a simple threefold singularity, and let $M'$
be a non-free MCM $A'$-module. We shall show that $M'$ does not admit a
$\g$-connection. Calculations in {\sc Singular}, mentioned in the
comments preceding the theorem, shows that the claim is true when $A'$
is of type $E_6$, $E_7$ and $E_8$, so we may assume that $A'$ is of
type $A_n$ or $D_n$. Moreover, we may assume that $M'$ is
indecomposable by lemma \ref{l:dirsum}.

Let $A = k[z_0,z_1,z_2,z_3]/(f)$ be the quasi-homogeneous $k$-algebra
that corresponds to the simple threefold singularity of type $A_n$ or
$D_n$, i.e. $\widehat A \cong A'$. By Yoshino \cite{yos90}, theorem
15.14, we may assume that $M' \cong \widehat M$ for some graded MCM
$A$-module $M$, and clearly $M$ is non-free and indecomposable. Since
$\g$ is generated by homogeneous derivations, it is a graded
Lie-Rinehart algebra. By lemma \ref{l:conn-is} and lemma
\ref{l:gr-conn}, it is therefore enough to show that if $M$ is a graded
indecomposable non-free MCM $A$-module, then $M$ does not admit a
graded $\g$-connection.

Assume that $\nabla: \g \to \enm_k(M)$ is a graded $\g$-connection on
$M$ for some graded indecomposable non-free MCM $A$-module $M$, and let
$d_0: L_1 \to L_0$ be a graded presentation of $M$. For any homogeneous
derivation $D \in \g$ of degree $\omega$, it follows from lemma
\ref{l:connpres} that $\nabla_D$ is induced by an operator of the form
$D + \phi_D: L_0 \to L_0$, where $\phi_D \in \enm_A(L_0)_\omega$
satisfies $D(d_0) = d_0 \psi_D - \phi_D d_0$ for some $\psi_D \in
\enm_A(L_1)_\omega$. Moreover, for any graded relation $\sum a_i D_i =
0$ in $\g$, we must have $\sum a_i \phi_{D_i} = 0$ in $\enm_A(M)$.

In the $A_n$ case, we choose coordinates such that $A =
k[x,y,z,w]/(f)$, where $f = x^2 + y^{n+1} + zw$ for some $n \ge 1$.
There are graded trivial derivations $D_1, D_2, D_3$ of $A$, all of
degree $0$, given by
\begin{equation*}
D_{1} = z \frac{\partial}{\partial z} - w \frac{\partial}{\partial w},
\quad D_{2} = w \frac{\partial}{\partial x} - 2x
\frac{\partial}{\partial z}, \quad D_{3} = z \frac{\partial}{\partial
x} - 2x \frac{\partial}{\partial w}
\end{equation*}
and there is a graded relation $2x D_1 + z D_2 - w D_3 = 0$. For all
graded indecomposable non-free MCM $A$-module $M$, this leads to a
contradiction, see appendix \ref{a:an} for details.

In the $D_n$ case, we choose coordinates such that $A =
k[x,y,z,w]/(f)$, where $f = x^2y + y^{n-1} + zw$ for some $n \ge 4$.
There are graded trivial derivations $D_1, D_2, D_3$ of $A$, with $D_1$
of degree $0$ and $D_2, D_3$ of degree $1$, given by
\begin{equation*}
D_{1} = z \frac{\partial}{\partial z} - w \frac{\partial}{\partial w},
\quad D_{2} = w \frac{\partial}{\partial x} - 2xy
\frac{\partial}{\partial z}, \quad D_{3} = z \frac{\partial}{\partial
x} - 2xy \frac{\partial}{\partial w}
\end{equation*}
and there is a graded relation $2xy D_1 + z D_2 - w D_3 = 0$. For
almost all graded indecomposable non-free MCM $A$-modules $M$, this
leads to a contradiction, see appendix \ref{a:an} for details.

However, in some exceptional cases, we shall use the graded trivial
derivations $D_4, D_5$ of $A$, with degree $n-3$, given by
\begin{equation*}
D_{4} = w \frac{\partial}{\partial y} - \beta \frac{\partial}{\partial
z}, \quad D_{5} = z \frac{\partial}{\partial y} - \beta
\frac{\partial}{\partial w}
\end{equation*}
where $\beta = x^2 + (n-1)y^{n-2}$, and the graded relation $\beta D_1
+ z D_4 - w D_5 = 0$. For the remaining graded indecomposable non-free
MCM $A$-modules $M$, this leads to a contradiction, see appendix
\ref{a:an} for details.
\end{proof}

The classification of non-Gorenstein singularities of finite CM
representation type is not known in dimension $d \ge 3$. However,
partial results are given in Eisenbud, Herzog \cite{ei-he88},
Auslander, Reiten \cite{au-re89} and Yoshino \cite{yos90}.

In Auslander, Reiten \cite{au-re89}, it was shown that there is only
one quotient singularity of dimension $d \ge 3$ with finite CM
representation type, the cyclic threefold quotient singularity of type
$\frac{1}{2} (1,1,1)$. Its complete local ring $A = S^G$ is the
invariant ring of the action of the group $G = \z_2$ on $S =
k[[x_1,x_2,x_3]]$ given by $\sigma x_i = -x_i$ for $i=1,2,3$, where
$\sigma \in G$ is the non-trivial element. There are exactly two
non-free indecomposable MCM $A$-modules $M_1$ and $M_2$. The module
$M_1$ has rank one, and is induced by the non-trivial representation of
$G$ of dimension one, see Auslander, Reiten \cite{au-re89}. By
proposition \ref{p:conn-quot}, it follows that $M_1$ admits an
integrable connection. One may show that $M_1$ is the canonical module
of $A$.

It is also known that the threefold scroll of type $(2,1)$, with
complete local ring $A = k[[x,y,z,u,v]]/(xz - y^2, xv - yu, yv - zu)$,
has finite CM representation type, see Auslander, Reiten
\cite{au-re89}. There are four non-free indecomposable MCM $A$-modules,
and none of these admit connections. In particular, the canonical
module $\omega_A$ does not admit a connection.

To the best of our knowledge, no other examples of singularities of
finite CM representation type are known in dimension $d \ge 3$.

\section{The canonical module}

Let $k$ be an algebraically closed field of characteristic $0$, and
let $A$ be a complete local CM $k$-algebra with residue field $k$
that is an isolated singularity. Then $A$ has a canonical module
$\omega_A$, an MCM $A$-module of finite injective dimension and with
type $r(\omega_A) = \dim_k \ext^d_A(k,\omega_A) = 1$, where $d =
\dim(A)$, see Bruns, Herzog \cite{br-he93}.

Since $A$ is CM, $\omega_A$ is also a dualizing module, see Bruns,
Herzog \cite{br-he93}, theorem 3.3.10. In particular,
$\ext^1_A(\omega_A, \omega_A) = 0$, and it follows from the
obstruction theory for connections that $\omega_A$ admits a
$k$-linear connection, see Eriksen, Gustavsen \cite{er-gu06}. If $A$
is Gorenstein, then $\omega_A$ is free, hence $\omega_A$ admits an
integrable connection. However, it turns out that this is not true
in general, and we ask the following question: \emph{When does the
canonical module admit a connection?}

When $A$ has dimension $d \ge 2$, it follows from our assumptions that
$A$ is normal. In this case, we say that $A$ is \emph{$\q$-Gorenstein}
if
    \[ \omega_A^{[n]} = ( \omega_A \otimes_A \dots \otimes_A
    \omega_A)^{\vee\vee} \cong A \]
for some integer $n \ge 1$. We remark that $A$ is $\q$-Gorenstein if
and only if $A$ is of the form $A = S^G$, where $S$ is the complete
local Gorenstein $k$-algebra with residue field $k$ and $G$ is a
finite subgroup of $\aut_k(S)$ that acts freely outside the closed
point of $\spec(S)$.

\begin{thm} \label{t:qgor-conn}
Let $A$ be the complete local ring of a singularity of dimension $d \ge
1$. If $A$ is a Gorenstein curve or a $\q$-Gorenstein singularity of
dimension $d \ge 2$, then the canonical module $\omega_A$ admits an
integrable connection.
\end{thm}
\begin{proof}
If $A$ is a Gorenstein curve, then $\omega_A \cong A$ and the result
is trivial. We may therefore assume that $A = S^G$ for a Gorenstein
singularity $S$ and a finite subgroup $G$ of $\aut_k(S)$ that acts
freely outside the closed point of $\spec(S)$. This implies that
there is a character $\chi$ of $G$ such that $\omega_A$ is
isomorphic to the semi-invariants $S^\chi = \{ a \in S: g a =
\chi(g) a \text{ for all } g \in G \}$, see Hini\v{c} \cite{hi76}.
Since $A$ is an isolated singularity of dimension at least two, the
canonical map $\der_k(S)^G \to \der_k(S^G)$ is an isomorphism, see
Schlessinger \cite{sch71}. Hence it follows from the proof of
proposition \ref{p:conn-quot} that $\omega_A$ admits an integrable
connection.
\end{proof}

We remark that it is \emph{not true} that the canonical module
$\omega_A$ admits an integrable connection when $A$ is a quotient of
the form $A = S^G$, where $S$ is a Gorenstein curve singularity and
$G$ is a finite subgroup of $\aut_k(S)$ that acts freely outside the
closed point of $\spec(S)$. In fact, consider the Gorenstein
monomial curve singularity $S = k[[t^3,t^5]]$, let $G = \z_2$, and
let the non-trivial element $\sigma \in \z_2$ act on $S$ as follows:
    \[ \sigma t^3 = -t^3, \quad \sigma t^5 = -t^5 \]
Then $A = S^G = k[[t^6,t^8,t^{10}]] \cong k[[t^3,t^4,t^5]]$, and we
see that $A$ is non-Gorenstein. This implies that $\omega_A$ does
not admit connections, as the following theorem shows:

\begin{thm} \label{t:canconn-mon}
Let $A$ be the complete local ring of a monomial curve singularity.
Then the canonical module $\omega_A$ admits a connection if and only if
$A$ is Gorenstein.
\end{thm}
\begin{proof}
The canonical module is characterized by its Hilbert function, and
we see from the characterization in Bruns, Herzog \cite{br-he93},
theorem 4.4.6 that the canonical module $\omega_A \cong
k[[\Lambda]]$, where $\Lambda = \Gamma \cup (\Gamma + \Delta)$
(using the notation of subsection \ref{ss:dim1}). Since $g \not \in
\Lambda$, the result follows from proposition \ref{p:nconn-mon} and
theorem \ref{t:conn-mon}.
\end{proof}

We would like to find necessary conditions for $\omega_A$ to admit a
connection, and the sufficient condition given in theorem
\ref{t:qgor-conn} is a natural candidate. In dimension one, we have
seen that the Gorenstein condition is necessary for any monomial
curve and for the curves $D^s_n$ for $n \le 100$ and $E^s_6, E^s_7,
E^s_8$. In dimension $d \ge 2$, we have seen that the
$\q$-Gorenstein condition is necessary for all known examples of
singularities of finite CM representation type. We make the
following conjecture:

\begin{conj}
Let $A$ be the complete local ring of a singularity of dimension $d \ge
1$. Then the canonical module $\omega_A$ admits a connection if and
only if $A$ is a Gorenstein curve or a $\q$-Gorenstein singularity of
dimension $d \ge 2$.
\end{conj}

\appendix

\section{Calculations for $A_n$ and $D_n$ in dimension 3}
\label{a:an}

Let $A$ be the quasi-homogeneous $k$-algebra corresponding to a
simple threefold singularity of type $A_n$ or $D_n$, and let $\g
\subseteq \der_k(A)$ denote the graded submodule generated by the
trivial derivations. In the proof of theorem \ref{t:simple-dim3}, we
used the fact that a graded non-free indecomposable MCM $A$-module
$M$ does not admit a graded $\g$-connection. There are complete
lists of such modules, and for each module $M$ in these lists, the
assumption that there exists a graded $\g$-connection on $M$ leads
to a contradiction. In this appendix, we shall give the full details
of the calculations that lead to these contradictions. We use the
notation of the proof of theorem \ref{t:simple-dim3}.

\subsection{The $A_n$ case}

Let $n \ge 1$ be an integer, and consider the quasi-homogeneous
$k$-algebra $A = k[x,y,z,w]/(f)$, where $f = x^2 + y^{n+1} + zw$ and
the grading of $A$ is given by
    \[ (\deg x, \deg y, \deg z, \deg w) = (n+1,2,n+1,n+1). \]
We consider the set of isomorphism classes of indecomposable graded
non-free MCM $A$-modules (up to degree shifting).

When $n$ is even, we may take the modules $\{ M_l: 1 \le l \le p\}$ as
representatives, where $p = \frac{n}{2}$ and the module $M_l$ has
presentation matrix
\[
\textsf M_l =  \begin{pmatrix} z & 0 & -x & -y^{n+1-l} \\ 0 & z &
-y^{l} & x \\ x & y^{n+1-l} & w & 0 \\ y^{l} & -x & 0 & w \end{pmatrix}
\]
When $n$ is odd, we may take the modules $\{ M_l: 1 \le l \le p-1
\}$ and $\{ N_-, N_+ \}$ as representatives, where $p =
\frac{n+1}{2}$, the module $M_l$ has presentation matrix as above,
and the modules $N_-$ and $N_+$ have presentation matrices
\begin{align*}
\textsf N_-  &  = \begin{pmatrix} z & -(x+iy^p) \\ x-iy^p & w
\end{pmatrix} &
\textsf N_+  &  = \begin{pmatrix} z & -(x-iy^p) \\ x+iy^p & w
\end{pmatrix}
\end{align*}
In each case, we view the presentation matrix as the matrix of a graded
$A$-linear map $d_0: L_1 \to L_0$ with respect to some chosen
homogeneous bases $\{ e_i \}$ for $L_0$ and $\{ f_i \}$ for $L_1$, with
degrees given by
\begin{align*}
( \deg e_1, \deg e_2, \deg e_3, \deg e_4 ) & = (0,n+1-2l,0,n+1-2l) \\
( \deg f_1, \deg f_2, \deg f_3, \deg f_4 ) & =(n+1,2n+2-2l,n+1,
2n+2-2l) \\
\intertext{for the modules $M_l$, and}
( \deg e_1, \deg e_2 ) & = (0,0) \\
( \deg f_1, \deg f_2 ) & = (n+1,n+1)
\end{align*}
for the modules $N_-$ and $N_+$.

\emph{The module $M_l$ for $n$ even}. Let $P$ and $C$ be matrices
corresponding to graded endomorphisms of $L_0$ and $L_1$ of degree $0$.
This implies that $\deg p_{ij} = \deg e_j - \deg e_i$ and $\deg c_{ij}
= \deg f_j - \deg f_i$, so $P$ and $C$ must be of the form
\begin{equation*}
P = \begin{pmatrix} p_{11} & 0 & p_{13} & 0 \\ 0 & p_{22} & 0 & p_{23}
\\ p_{31} & 0 & p_{33} & 0 \\ 0 & p_{42} & 0 & p_{44} \end{pmatrix}, \;
C = \begin{pmatrix} c_{11} & 0 & c_{13} & 0 \\ 0 & c_{22} & 0 & c_{23}
\\ c_{31} & 0 & c_{33} & 0 \\ 0 & c_{42} & 0 & c_{44}
\end{pmatrix}
\end{equation*}
with $p_{ij}, c_{ij} \in k$ for $1 \le i,j \le 4$. For $s = 1,2,3$, we
must solve the equation $( D_s(\textsf M_l) ) = \textsf M_l C_s - P_s
\textsf M_l$, where $D_s$ is the derivation of $A$ mentioned in the
proof of theorem \ref{t:simple-dim3}, and $P_s$ and $C_s$ have the
above form. This gives $P_s = P_s^0 + a_s \Psi_0$ with $a_s \in k$ for
$s=1,2,3$, where $\Psi_0 = I_4$ and
\begin{equation*}
P_1^0 = \begin{pmatrix} 0 & 0 & 0 & 0 \\ 0 & 0 & 0 & 0 \\ 0 & 0 & 1 &
0 \\ 0 & 0 & 0 & 1 \end{pmatrix}, P_2^0 = \begin{pmatrix} 0 & 0 & -1 &
0 \\ 0 & 0 & 0 & 1 \\ 0 & 0 & 0 & 0 \\ 0 & 0 & 0 & 0 \end{pmatrix},
P_3^0 = \begin{pmatrix} 0 & 0 & 0 & 0 \\ 0 & 0 & 0 & 0 \\
1 & 0 & 0 & 0 \\ 0 & -1 & 0 & 0 \end{pmatrix}
\end{equation*}
Let $\alpha = 2x a_{1} + z a_2 - w a_3$. Using the relation $2x D_1
+ z D_2 - w D_3 = 0$, we obtain the equation
    \[ 2x P_1 + z P_2 - w P_3 = \begin{pmatrix} \alpha & 0 & -z & 0 \\ 0 &
    \alpha & 0 & z \\ -w & 0 & 2x + \alpha & 0 \\ 0 & w & 0 & 2x +
    \alpha \end{pmatrix} = 0  \]
in $\enm_A(M_l)$. By inspection, we see that this is a contradiction
for $1 \le l \le p$.

\emph{The module $M_l$ for $n$ odd}. Let $P$ and $C$ be matrices
corresponding to graded endomorphisms of $L_0$ and $L_1$ of degree $0$.
This implies that $\deg p_{ij} = \deg e_j - \deg e_i$ and $\deg c_{ij}
= \deg f_j - \deg f_i$, so $P$ and $C$ must be of the form
\begin{equation*}
P = \begin{pmatrix} p_{11} & p_{12} y^{p-l} & p_{13} & p_{14} y^{p-l}
\\ 0 & p_{22} & 0 & p_{23} \\ p_{31} & p_{32} y^{p-l} & p_{33} & p_{34}
y^{p-l} \\ 0 & p_{42} & 0 & p_{44} \end{pmatrix}, \; C =
\begin{pmatrix} c_{11} & c_{12} y^{p-l} & c_{13} & c_{14} y^{p-l} \\ 0
& c_{22} & 0 & c_{23} \\ c_{31} & c_{32} y^{p-l} & c_{33} & c_{34}
y^{p-l} \\ 0 & c_{42} & 0 & c_{44}
\end{pmatrix}
\end{equation*}
with $p_{ij}, c_{ij} \in k$ for $1 \le i,j \le 4$. For $s = 1,2,3$, we
must solve the equation $( D_s(\textsf M_l) ) = \textsf M_l C_s - P_s
\textsf M_l$, where $D_s$ is the derivation of $A$ mentioned in the
proof of theorem \ref{t:simple-dim3}, and $P_s$ and $C_s$ have the
above form. This gives $P_s = P_s^0 + a_s \Psi_0$ with $a_s \in k$ for
$s=1,2,3$, where $\Psi_0 = I_4$ and
\begin{equation*}
P_1^0 = \begin{pmatrix} 0 & 0 & 0 & 0 \\ 0 & 0 & 0 & 0 \\ 0 & 0 & 1 &
0 \\ 0 & 0 & 0 & 1 \end{pmatrix}, P_2^0 = \begin{pmatrix} 0 & 0 & -1 &
0 \\ 0 & 0 & 0 & 1 \\ 0 & 0 & 0 & 0 \\ 0 & 0 & 0 & 0 \end{pmatrix},
P_3^0 = \begin{pmatrix} 0 & 0 & 0 & 0 \\ 0 & 0 & 0 & 0 \\
1 & 0 & 0 & 0 \\ 0 & -1 & 0 & 0 \end{pmatrix}
\end{equation*}
As in the case of the modules $M_l$ for $n$ even, the equation $2x
P_1 + z P_2 - w P_3 = 0$ in $\enm_A(M_l)$ leads to a contradiction
for $1 \le l \le p-1$.

\emph{The module $N_-$ for $n$ odd}. Let $P$ and $C$ be matrices
corresponding to graded endomorphisms of $L_0$ and $L_1$ of degree $0$.
This implies that $\deg p_{ij} = \deg e_j - \deg e_i$ and $\deg c_{ij}
= \deg f_j - \deg f_i$, so $P$ and $C$ must be of the form
\begin{equation*}
P = \begin{pmatrix} p_{11} & p_{12} \\ p_{21} & p_{22} \end{pmatrix},
\; C = \begin{pmatrix} c_{11} & c_{12} \\ c_{21} & c_{22} \end{pmatrix}
\end{equation*}
with $p_{ij}, c_{ij} \in k$ for $1 \le i,j \le 2$. For $s = 1,2,3$, we
must solve the equation $( D_s(\textsf N_-) ) = \textsf N_- C_s - P_s
\textsf N_-$, where $D_s$ is the derivation of $A$ mentioned in the
proof of theorem \ref{t:simple-dim3}, and $P_s$ and $C_s$ have the
above form. This gives $P_s = P_s^0 + a_s \Psi_0$ with $a_s \in k$ for
$s=1,2,3$, where $\Psi_0 = I_2$ and
\begin{equation*}
P_1^0 = \begin{pmatrix} 0 & 0 \\ 0 & -1 \end{pmatrix}, P_2^0 =
\begin{pmatrix} 0 & -1 \\ 0 & 0 \end{pmatrix},
P_3^0 = \begin{pmatrix} 0 & 0 \\ 1 & 0 \end{pmatrix}
\end{equation*}
Let $\alpha = 2x a_{1} + z a_2 - w a_3$. Using the relation $2x D_1
+ z D_2 - w D_3 = 0$, we obtain the equation
    \[ 2x P_1 + z P_2 - w P_3 = \begin{pmatrix} \alpha & z \\ w &
    -2x + \alpha \end{pmatrix} = 0 \]
in $\enm_A(N_-)$. By inspection, we see that this is a
contradiction.

\emph{The module $N_+$ for $n$ odd}. We see that $N_+ \cong N_-^\vee$
using lemma \ref{l:mfid}. It follows from lemma \ref{l:dual} and the
computation above for $N_-$ that the module $N_+$ cannot admit a
$\g$-connection.

\subsection{The $D_n$ case}

Let $n \ge 4$ be an integer, and consider the quasi-homogeneous
$k$-algebra $A = k[x,y,z,w]/(f)$, where $f = x^2y + y^{n-1} + zw$ and
the grading of $A$ is given by
    \[ (\deg x, \deg y, \deg z, \deg w) = (n-2,2,n-1,n-1). \]
We consider the set of isomorphism classes of indecomposable graded
non-free MCM $A$-modules (up to degree shifting).

When $n$ is odd, we may take the modules $\{ M_l: 1 \le l \le p \}$,
$\{ N_l: 1 \le l \le p \}$, $\{ X_l: 1 \le l \le p+1 \}$, $\{ Y_l: 1
\le l \le p \}$ and $\{ B_1, B_2 \}$ as representatives, where $p =
\frac{n-3}{2}$ and the modules $M_l, N_l, X_l, Y_l, B_1, B_2$ have
presentation matrices
\begin{align*}
\textsf M_l & = \begin{pmatrix} z & 0 & -xy & -y^{n-1-l} \\ 0 &
z & -y^{l+1} & xy \\ x & y^{n-2-l} & w & 0 \\ y^l & -x & 0 & w
\end{pmatrix}  &
\textsf N_l & = \begin{pmatrix} z & 0 & -x & -y^{n-2-l} \\ 0 & z &
-y^{l} & x \\ xy & y^{n-1-l} & w & 0 \\ y^{l+1} & -xy & 0 & w
\end{pmatrix} \\
\textsf X_l & = \begin{pmatrix} z & 0 & -x & -y^{n-1-l} \\ 0 &
z & -y^{l} & xy \\ xy & y^{n-1-l} & w & 0 \\ y^l & -x & 0 & w
\end{pmatrix} &
\textsf Y_l & = \begin{pmatrix} z & 0 & -xy & -y^{n-1-l} \\ 0 &
z & -y^{l} & x \\ x & y^{n-1-l} & w & 0 \\ y^l & -xy & 0 & w
\end{pmatrix} \\
\textsf B_1  &  = \begin{pmatrix} z & -(x^2+y^{n-2}) \\ y & w
\end{pmatrix} &
\textsf B_2  &  = \begin{pmatrix} z & -y \\ x^2+y^{n-2} & w
\end{pmatrix}
\end{align*}
When $n$ is even, we may take the modules $\{ M_l: 1 \le l \le p-1
\}$, $\{ N_l: 1 \le l \le p-1 \}$, $\{ X_l: 1 \le l \le p \}$, $\{
Y_l: 1 \le l \le p \}$ and $\{ B_1, B_2, C_-, C_+, D_-, D_+ \}$ as
representatives, where $p = \frac{n-2}{2}$, the modules $M_l, N_l,
X_l, Y_l, B_1, B_2$ have presentation matrices as above, and the
modules $C_-, C_+, D_-, D_+$ have presentation matrices
\begin{align*}
\mathsf C_- & = \begin{pmatrix} z & -(x+iy^p) \\ y(x-iy^p) & w
\end{pmatrix} &
\mathsf C_+ & = \begin{pmatrix} z & -(x-iy^p) \\ y(x+iy^p) & w
\end{pmatrix} \\
\mathsf D_- & = \begin{pmatrix} z & -y(x+iy^p) \\ x-iy^p & w
\end{pmatrix} &
\mathsf D_+ & = \begin{pmatrix} z & -y(x-iy^p) \\ x+iy^p & w
\end{pmatrix}
\end{align*}
For each of these representatives, we view the presentation matrix as
the matrix of a graded $A$-linear map $d_0: L_1 \to L_0$ with respect
to some chosen homogeneous bases $\{ e_i \}$ for $L_0$ and $\{ f_i \}$
for $L_1$, with degrees given by
\begin{align*}
( \deg e_1, \deg e_2, \deg e_3, \deg e_4 ) & = (0,n-2-2l,1,n-1-2l) \\
( \deg f_1, \deg f_2, \deg f_3, \deg f_4 ) & =(n-1,2n-3-2l,n,
2n-2-2l) \\
\intertext{for the module $M_l$,}
( \deg e_1, \deg e_2, \deg e_3, \deg e_4 ) & = (0,n-2-2l,-1,n-3-2l) \\
( \deg f_1, \deg f_2, \deg f_3, \deg f_4 ) & =(n-1,2n-3-2l,n-2,
2n-4-2l) \\
\intertext{for the module $N_l$,}
( \deg e_1, \deg e_2, \deg e_3, \deg e_4 ) & = (0,n-2-2l,-1,n-1-2l) \\
( \deg f_1, \deg f_2, \deg f_3, \deg f_4 ) & =(n-1,2n-3-2l,n-2,
2n-2-2l) \\
\intertext{for the module $X_l$,}
( \deg e_1, \deg e_2, \deg e_3, \deg e_4 ) & = (0,n-2l,1,n-1-2l) \\
( \deg f_1, \deg f_2, \deg f_3, \deg f_4 ) & =(n-1,2n-1-2l,n,
2n-2-2l) \\
\intertext{for the module $Y_l$,}
( \deg e_1, \deg e_2 ) & = (0,n-3) \\
( \deg f_1, \deg f_2 ) & = (n-1,2n-4)
\intertext{for the module $B_1$,}
( \deg e_1, \deg e_2 ) & = (0,3-n) \\
( \deg f_1, \deg f_2 ) & = (n-1,2)
\intertext{for the module $B_2$,}
( \deg e_1, \deg e_2 ) & = (0,-1) \\
( \deg f_1, \deg f_2 ) & = (n-1,n-2)
\intertext{for the modules $C_-$ and $C_+$, and }
( \deg e_1, \deg e_2 ) & = (0,1) \\
( \deg f_1, \deg f_2 ) & = (n-1,n)
\intertext{for the modules $D_-$ and $D_+$.}
\end{align*}

\emph{The module $M_l$ for $n$ odd}. Let $P$ and $C$ be matrices
corresponding to graded endomorphisms of $L_0$ and $L_1$ of degree
$\omega$. Then $\deg p_{ij} = \deg e_j - \deg e_i + \omega$ and $\deg
c_{ij} = \deg f_j - \deg f_i + \omega$. In case $\omega = 0$, $P$ and
$C$ must be of the form
\begin{align*}
P & = \begin{pmatrix} p_{11} & 0 & 0 & p_{14} y^{p-l+1} \\ 0 & p_{22} &
p_{23} y^{l-p} & 0 \\ 0 & p_{32} y^{p-l} & p_{33} & 0 \\ 0 & 0 & 0 &
p_{44} \end{pmatrix} \\
C & = \begin{pmatrix} c_{11} & 0 & 0 & c_{14} y^{p-l+1} \\ 0 & c_{22} &
c_{23} y^{l-p} & 0 \\ 0 & c_{32} y^{p-l} & c_{33} & 0 \\ 0 & 0 & 0 &
c_{44} \end{pmatrix}
\end{align*}
with $p_{ij}, c_{ij} \in k$ for $1 \le i,j \le 4$, and $p_{23} = c_{23}
= 0$ if $l \neq p$. In case $\omega = 1$, $P$ and $C$ must be of the
form
\begin{align*}
P & = \begin{pmatrix} 0 & p_{12} y^{p-l+1} & p_{13} y & p_{14} x \\
p_{21} & 0 & 0 & p_{24} y \\ p_{31} & 0 & 0 & p_{34} y^{p-l+1} \\
0 & p_{42} & p_{43} & 0 \end{pmatrix} \\
C & = \begin{pmatrix} 0 & c_{12} y^{p-l+1} & c_{13} y & c_{14} x \\
c_{21} & 0 & 0 & c_{24} y \\ c_{31} & 0 & 0 & c_{34} y^{p-l+1} \\
0 & c_{42} & c_{43} & 0 \end{pmatrix}
\end{align*}
with $p_{ij}, c_{ij} \in k$ for $1 \le i,j \le 4$, and $p_{21} = c_{21}
= p_{43} = c_{43} = 0$ if $l \neq p$ and $p_{14} = c_{14} = 0$ if $l
\neq 1$. We must solve the equation $( D_s(\textsf M_l) ) = \textsf M_l
C_s - P_s \textsf M_l$ for $s=1,2,3$, where $D_s$ is the derivation of
$A$ mentioned in the proof of theorem \ref{t:simple-dim3}, and $P_s$
and $C_s$ have the above form (with $\omega = 0$ for $s = 1$ and
$\omega = 1$ for $s = 2,3$). This gives $P_1 = P_1^0 + a_1 \Phi_0$ with
$a_1 \in k$, where $\Phi_0 = I_4$ and
\begin{equation*}
P_1^0 = \begin{pmatrix} -1 & 0 & 0 & 0 \\ 0 & -1 & 0 & 0 \\ 0 & 0 & 0 &
0 \\ 0 & 0 & 0 & 0 \end{pmatrix}
\end{equation*}
Moreover, $P_s = P_s^0 + a_s \Phi_1$ with $a_s \in k$ for $s = 2,3$,
and $a_2 = a_3 = 0$ if $l \neq p$, where
\begin{equation*}
\Phi_1 = \begin{pmatrix} 0 & y & 0 & 0 \\ -1 & 0 & 0 & 0 \\ 0 & 0 & 0 &
-y \\0 & 0 & 1 & 0 \end{pmatrix}
\end{equation*}
and
\begin{equation*}
P_2^0 = \begin{pmatrix} 0 & 0 & y & 0 \\ 0 & 0 & 0 & -y \\ 0 & 0 & 0 &
0 \\ 0 & 0 & 0 & 0 \end{pmatrix}, P_3^0 = \begin{pmatrix} 0 & 0 & 0 & 0
\\ 0 & 0 & 0 & 0 \\ -1 & 0 & 0 & 0 \\ 0 & 1 & 0 & 0 \end{pmatrix}
\end{equation*}
The relation $2xy D_1 + z D_2 - w D_3 = 0$ gives the equation $2xy
P_1 + z P_2 - w P_3 = 0$ in $\enm_A(M_l)$, i.e.
\begin{equation*}
\begin{pmatrix} 2(a_1-1)xy & a_2 yz - a_3 yw & yz & 0 \\ -a_2 z + a_3 w
& 2(a_1-1) xy & 0 & -yz \\ w & 0 & 2a_1 xy & -a_2 yz + a_3 yw \\ 0 &
-w & a_2 y - a_3 w & 2 a_1 xy \end{pmatrix} = 0
\end{equation*}
By inspection, we see that this is a contradiction for $1 \le l \le
p$.

\emph{The module $N_l$ for $n$ odd}. We see that $N_l \cong M_l^\vee$
for $1 \le l \le p$ using lemma \ref{l:mfid}. It follows from lemma
\ref{l:dual} and the computation above for $M_l$ that the module $N_l$
cannot admit a $\g$-connection for $1 \le l \le p$.

\emph{The module $X_l$ for $n$ odd}. We see that $X_l \cong Y_l^\vee$
for $1 \le l \le p$ using lemma \ref{l:mfid}. It follows from lemma
\ref{l:dual} and the computation below for $Y_l$ that the module $X_l$
cannot admit a $\g$-connection for $1 \le l \le p$. Since we include
the case $l=p+1$ in the calculations below for $Y_l$, it also follows
that the module $X_{p+1} \cong Y_{p+1}$ cannot admit a $\g$-connection.

\emph{The module $Y_l$ for $n$ odd}. In this case, we include the
module $Y_l$ for $1 \le l \le p+1$ for reasons mentioned above. Let $P$
and $C$ be matrices corresponding to graded endomorphisms of $L_0$ and
$L_1$ of degree $\omega$. Then $\deg p_{ij} = \deg e_j - \deg e_i +
\omega$ and $\deg c_{ij} = \deg f_j - \deg f_i + \omega$. In case
$\omega = 0$, $P$ and $C$ must be of the form
\begin{align*}
P & = \begin{pmatrix} p_{11} & p_{12} x & 0 & p_{14} y^{p-l+1} \\ 0 &
p_{22} & p_{23} & 0 \\ 0 & p_{32} y^{p+1-l} & p_{33} & 0 \\ p_{41} & 0
& 0 & p_{44} \end{pmatrix} \\
C & = \begin{pmatrix} c_{11} & c_{12} x & 0 & c_{14} y^{p-l+1} \\ 0 &
c_{22} & c_{23} & 0 \\ 0 & c_{32} y^{p+1-l} & c_{33} & 0 \\ c_{41} & 0
& 0 & c_{44} \end{pmatrix}
\end{align*}
with $p_{ij}, c_{ij} \in k$ for $1 \le i,j \le 4$, and $p_{23} = c_{23}
= p_{41} = c_{41} = 0$ if $l \neq p+1$, $p_{12} = c_{12} = 0$ if $l
\neq 1$. In case $\omega = 1$, $P$ and $C$ must be of the form
\begin{align*}
P & = \begin{pmatrix} 0 & p_{12} y^{p-l+2} + p_{12}' z + p_{12}'' w &
p_{13} y & p_{14} x \\ p_{21} & 0 & 0 & p_{24} \\ p_{31} & p_{32} x & 0
& p_{34} y^{p+1-l} \\ 0 & p_{42} y & p_{43} y^{l-p} & 0 \end{pmatrix} \\
C & = \begin{pmatrix} 0 & c_{12} y^{p-l+2} + c_{12}' z + c_{12}'' w &
c_{13} y & c_{14} x \\ c_{21} & 0 & 0 & c_{24} \\ c_{31} & c_{32} x & 0
& c_{34} y^{p+1-l} \\ 0 & c_{42} y & c_{43} y^{l-p} & 0 \end{pmatrix}
\end{align*}
with $p_{ij}, p_{ij}', p_{ij}'', c_{ij}, c_{ij}', c_{ij}'' \in k$ for
$1 \le i,j \le 4$, and $p_{43} = c_{43} = 0$ if $l < p$, $p_{12}' =
p_{12}'' = c_{12}' = c_{12}'' = p_{14} = c_{14} = p_{32} = c_{32} = 0$
if $l \neq 1$, $p_{21} = c_{21} = 0$ if $l \neq p+1$. For $s = 1,2,3$,
we must solve the equation $( D_s(\textsf Y_l) ) = \textsf Y_l C_s -
P_s \textsf Y_l$, where $D_s$ is the derivation of $A$ mentioned in the
proof of theorem \ref{t:simple-dim3}, and $P_s$ and $C_s$ have the
above form (with $\omega = 0$ for $s = 1$ and $\omega = 1$ for $s =
2,3$). This gives $P_1 = P_1^0 + a_1 \Phi_0$ with $a_1 \in k$, where
$\Phi_0 = I_4$ and
\begin{equation*}
P_1^0 = \begin{pmatrix} 0 & 0 & 0 & 0 \\ 0 & 0 & 0 & 0 \\ 0 & 0 & 1 & 0
\\ 0 & 0 & 0 & 1 \end{pmatrix}
\end{equation*}
Moreover, $P_s = P_s^0 + a_s \Psi_1 + a_s' \Psi_1' + a_s'' \Psi_1''$
for $s = 2,3$, and $a_s = a_s' = 0$ if $l \neq 1$, $a_s'' = 0$ if $l
\neq p+1$, where
\begin{equation*}
\Psi_1 = \begin{pmatrix} 0 & -w & -y & x \\ 0 & 0 & 0 & 0 \\ 0 & 0 & 0
& 0 \\ 0 & 0 & 0 & 0 \end{pmatrix}, \Psi_1' = \begin{pmatrix} 0 & z & 0
& 0 \\ 0 & 0 & 0 & 0 \\ 0 & x & 0 & 0 \\ 0 & y & 0 & 0 \end{pmatrix},
\Psi_1'' = \begin{pmatrix} 0 & y & 0 & 0 \\ -1 & 0 & 0 & 0 \\ 0 & 0 & 0
& -1 \\ 0 & 0 & y & 0 \end{pmatrix}
\end{equation*}
and
\begin{equation*}
P_2^0 = \begin{pmatrix} 0 & 0 & y & 0 \\ 0 & 0 & 0 & -1 \\ 0 & 0 & 0 &
0 \\ 0 & 0 & 0 & 0 \end{pmatrix}, P_3^0 = \begin{pmatrix} 0 & 0 & 0 & 0
\\ 0 & 0 & 0 & 0 \\ -1 & 0 & 0 & 0 \\ 0 & y & 0 & 0 \end{pmatrix}
\end{equation*}
The relation $2xy D_1 + z D_2 - w D_3 = 0$ give the equation $2xy
P_1 + z P_2 - w P_3 = 0$ in $\enm_A(Y_l)$, i.e.
\begin{equation*}
\begin{pmatrix} 2 a_1 xy & -bw+ b'z + b''y & yz - by & bx \\ -b'' & 2
a_1 xy & 0 & -1 \\ w & b'x & 2(a_1+1) xy & -b'' \\ 0 & -yw + b'y &
ab''y & 2 (a_1+1) xy \end{pmatrix} = 0,
\end{equation*}
where $b = a_2 z - a_3 w$, $b' = a_2' z - a_3' w$ and $b'' = a_2'' z
- a_3'' w$. By inspection, we see that this is a contradiction for
$1 \le l \le p+1$.

\emph{The module $B_1$ for $n$ odd}. Let $P$ and $C$ be matrices
corresponding to graded endomorphisms of $L_0$ and $L_1$ of degree
$\omega$. Then $\deg p_{ij} = \deg e_j - \deg e_i + \omega$ and $\deg
c_{ij} = \deg f_j - \deg f_i + \omega$. In case $\omega = 0$, $P$ and
$C$ must be of the form
\begin{equation*}
P = \begin{pmatrix} p_{11} & p_{12}y^{p} \\ 0 & p_{22} \end{pmatrix},
\; C = \begin{pmatrix} c_{11} & c_{12}y^p \\ 0 & c_{22} \end{pmatrix}
\end{equation*}
with $p_{ij}, c_{ij} \in k$ for $1 \le i,j \le 2$. In case $\omega =
n-3$, $P$ and $C$ must be of the form
\begin{align*}
P & = \begin{pmatrix} p_{11} y^p & p_{12} y^{n-3} + p_{12}' y^{p-1} w +
p_{12}'' y^{p-1} z \\ p_{21} & p_{22} y^p \end{pmatrix} \\
C & = \begin{pmatrix} c_{11} y^p & c_{12} y^{n-3} + c_{12}' y^{p-1} w +
c_{12}'' y^{p-1} z \\ c_{21} & c_{22} y^p \end{pmatrix}
\end{align*}
with $p_{ij}, p_{ij}', p_{ij}'', c_{ij}, c_{ij}', c_{ij}'' \in k$ for
$1 \le i,j \le 2$. We must solve the equation $( D_s(\textsf B_1) ) =
\textsf B_1 C_s - P_s \textsf B_1$ for $s = 1,4,5$, where $D_s$ is the
derivation of $A$ mentioned in the proof of theorem
\ref{t:simple-dim3}, and $P_s$ and $C_s$ have the above form (with
$\omega = 0$ for $s = 1$ and $\omega = n-3$ for $s = 4,5$). This gives
$P_1 = P_1^0 + a_1 \Phi_0$ with $a_1 \in k$, where $\Phi_0 = I_2$ and
\begin{equation*}
P_1^0 = \begin{pmatrix} 1 & 0 \\ 0 & 0 \end{pmatrix}
\end{equation*}
Moreover, $P_s = P_s^0 + a_s \Phi_{n-3}$ with $a_s \in k$ for $s =
4,5$, where we have $\Phi_{n-3} = y^p I_2$, $\gamma = (n-2) y^{n-3}$,
and
\begin{equation*}
P_4^0 = \begin{pmatrix} 0 & -\gamma \\ 0 & 0 \end{pmatrix}, \; P_5^0 =
\begin{pmatrix} 0 & 0 \\ 1 & 0 \end{pmatrix}
\end{equation*}
Let $\beta = x^2 + (n-1) y^{n-2}$. The relation $\beta D_1 + zD_4 -
wD_5 = 0$ gives the equation $\beta P_1 + z P_4 - w P_5 = 0$ in
$\enm_A(B_1)$, i.e.
\begin{equation*}
\begin{pmatrix} \beta + \delta & -z \gamma \\ -w & \delta
\end{pmatrix} = 0,
\end{equation*}
where $\delta = a_1 \beta + a_4 y^p z - a_5 y^p w$. By inspection,
we see that this is a contradiction.

\emph{The module $B_2$ for $n$ odd}. We see that $B_2 \cong B_1^\vee$
using lemma \ref{l:mfid}. It follows from lemma \ref{l:dual} and the
computation above for $B_1$ that the module $B_2$ cannot admit a
$\g$-connection.

\emph{The module $M_l$ for $n$ even}. Let $P$ and $C$ be matrices
corresponding to graded endomorphisms of $L_0$ and $L_1$ of degree
$\omega$. Then $\deg p_{ij} = \deg e_j - \deg e_i + \omega$ and $\deg
c_{ij} = \deg f_j - \deg f_i + \omega$. In case $\omega = 0$, $P$ and
$C$ must be of the form
\begin{equation*}
P = \begin{pmatrix} p_{11} & p_{12} y^{p-l} & 0 & 0 \\ 0 & p_{22} & 0 &
0 \\ 0 & 0 & p_{33} & p_{34} y^{p-l} \\ 0 & 0 & 0 & p_{44}
\end{pmatrix},
C = \begin{pmatrix} c_{11} & c_{12} y^{p-l} & 0 & 0 \\ 0 & c_{22} & 0 &
0 \\ 0 & 0 & c_{33} & c_{34} y^{p-l} \\ 0 & 0 & 0 & c_{44}
\end{pmatrix}
\end{equation*}
with $p_{ij}, c_{ij} \in k$ for $1 \le i,j \le 4$. In case $\omega =
1$, $P$ and $C$ must be of the form
\begin{align*}
P & = \begin{pmatrix} 0 & 0 & p_{13} y & p_{14} y^{p+1-l} + p_{14}' x
\\ 0 & 0 & p_{23} & p_{24} y \\ p_{31} & p_{32} y^{p-l} & 0 & 0 \\
0 & p_{42} & 0 & 0 \end{pmatrix} \\
C & = \begin{pmatrix} 0 & 0 & c_{13} y & c_{14} y^{p+1-l} + c_{14}' x
\\ 0 & 0 & c_{23} & c_{24} y \\ c_{31} & c_{32} y^{p-l} & 0 & 0 \\
0 & c_{42} & 0 & 0 \end{pmatrix}
\end{align*}
with $p_{ij}, p_{ij}', c_{ij}, c_{ij}' \in k$ for $1 \le i,j \le 4$,
and $p_{14}' = c_{14}' = 0$ if $l \neq 1$, $p_{23} = c_{23} = 0$ if $l
\neq p-1$. For $s = 1,2,3$, we must solve the equation $( D_s(\textsf
M_l) ) = \textsf M_l C_s - P_s \textsf M_l$, where $D_s$ is the
derivation of $A$ mentioned in the proof of theorem
\ref{t:simple-dim3}, and $P_s$ and $C_s$ have the above form (with
$\omega = 0$ for $s = 1$ and $\omega = 1$ for $s = 2,3$). This gives
$P_1 = P_1^0 + a_1 \Phi_0$ with $a_1 \in k$, where $\Phi_0 = I_4$ and
\begin{equation*}
P_1^0 = \begin{pmatrix} -1 & 0 & 0 & 0 \\ 0 & -1 & 0 & 0 \\ 0 & 0 & 0 &
0 \\ 0 & 0 & 0 & 0 \end{pmatrix}
\end{equation*}
Moreover, $P_s = P_s^0$ for $s = 2,3$, where
\begin{equation*}
P_2^0 = \begin{pmatrix} 0 & 0 & y & 0 \\ 0 & 0 & 0 & -y \\ 0 & 0 & 0 &
0 \\ 0 & 0 & 0 & 0 \end{pmatrix}, P_3^0 = \begin{pmatrix} 0 & 0 & 0 & 0
\\ 0 & 0 & 0 & 0 \\ -1 & 0 & 0 & 0 \\ 0 & 1 & 0 & 0 \end{pmatrix}
\end{equation*}
The relation $2xy D_1 + z D_2 - w D_3 = 0$ gives the equation $2xy
P_1 + z P_2 - w P_3 = 0$ in $\enm_A(M_l)$, i.e.
\begin{equation*}
\begin{pmatrix} 2(a_1-1)xy & 0 & yz & 0 \\ 0 & 2(a_1-1) xy & 0 & -yz \\
w & 0 & 2a_1 xy & 0 \\ 0 & -w & 0 & 2 a_1 xy \end{pmatrix} = 0
\end{equation*}
By inspection, we see that this is a contradiction for $1 \le l \le
p-1$.

\emph{The module $N_l$ for $n$ even}. We see that $N_l \cong M_l^\vee$
for $1 \le l \le p-1$ using lemma \ref{l:mfid}. It follows from lemma
\ref{l:dual} and the computation above for $M_l$ that the module $N_l$
cannot admit a $\g$-connection.

\emph{The module $X_l$ for $n$ even}. We see that $X_l \cong Y_l^\vee$
for $1 \le l \le p$ using lemma \ref{l:mfid}. It follows from lemma
\ref{l:dual} and the computation below for $Y_l$ that the module $X_l$
cannot admit a $\g$-connection.

\emph{The module $Y_l$ for $n$ even}. Let $P$ and $C$ be matrices
corresponding to graded endomorphisms of $L_0$ and $L_1$ of degree
$\omega$. Then $\deg p_{ij} = \deg e_j - \deg e_i + \omega$ and $\deg
c_{ij} = \deg f_j - \deg f_i + \omega$. In case $\omega = 0$, $P$ and
$C$ must be of the form
\begin{align*}
P & = \begin{pmatrix} p_{11} & p_{12} y^{p+1-l} + p_{12}' x & 0 & 0 \\
0 & p_{22} & 0 & 0 \\ 0 & 0 & p_{33} & p_{34} y^{p-l} \\ 0 &
0 & p_{43} & p_{44} \end{pmatrix} \\
C & = \begin{pmatrix} c_{11} & c_{12} y^{p+1-l} + c_{12}' x & 0 & 0 \\
0 & c_{22} & 0 & 0 \\ 0 & 0 & c_{33} & c_{34} y^{p-l} \\ 0 & 0 & c_{43}
& c_{44} \end{pmatrix}
\end{align*}
with $p_{ij}, p_{ij}', c_{ij}, c_{ij}' \in k$ for $1 \le i,j \le 4$,
and $p_{12}' = c_{12}'$ if $l \neq 1$, $p_{43} = c_{43} = 0$ if $l \neq
p$. In case $\omega = 1$, $P$ and $C$ must be of the form
\begin{align*}
P & = \begin{pmatrix} 0 & p_{12} z + p_{12}' w & p_{13} y + p_{13}' x &
p_{14} y^{p+1-l} + p_{14}' x \\ 0 & 0 & p_{23} & p_{24} \\ p_{31} &
p_{32} y^{p+1-l} + p_{32}' x & 0 & 0 \\ p_{41} & p_{42} y + p_{42}' x &
0 & 0 \end{pmatrix} \\
C & = \begin{pmatrix} 0 & c_{12} z + c_{12}' w & c_{13} y + c_{13}' x &
c_{14} y^{p+1-l} + c_{14}' x \\ 0 & 0 & c_{23} & c_{24} \\ c_{31} &
c_{32} y^{p+1-l} + c_{32}' x & 0 & 0 \\ c_{41} & c_{42} y + c_{42}' x &
0 & 0 \end{pmatrix}
\end{align*}
with $p_{ij}, p_{ij}', c_{ij}, c_{ij}' \in k$ for $1 \le i,j \le 4$,
and $p_{23} = c_{23} = p_{41} = c_{41} = 0$ if $l \neq p$, $p_{12} =
p_{12}' = c_{12} = c_{12}' = p_{14}' = c_{14}' = p_{32}' = c_{32}' = 0$
if $l \neq 1$, and furthermore $p_{13}' = c_{13}' = p_{42}' = c_{42}' =
0$ if $n \neq 4$. For $s = 1,2,3$, we must solve the equation $(
D_s(\textsf Y_l) ) = \textsf Y_l C_s - P_s \textsf Y_l$, where $D_s$ is
the derivation of $A$ mentioned in the proof of theorem
\ref{t:simple-dim3}, and $P_s$ and $C_s$ have the above form (with
$\omega = 0$ for $s = 1$ and $\omega = 1$ for $s = 2,3$). This gives
$P_1 = P_1^0 + a_1 \Phi_0$ with $a_1 \in k$, where $\Phi_0 = I_4$ and
\begin{equation*}
P_1^0 = \begin{pmatrix} -1 & 0 & 0 & 0 \\ 0 & -1 & 0 & 0 \\ 0 & 0 & 0 &
0 \\ 0 & 0 & 0 & 0 \end{pmatrix}
\end{equation*}
Moreover, $P_s = P_s^0 + a_s \Psi_1 + a_s' \Psi_1'$ with $a_s, a_s' \in
k$ for $s = 2,3$, and $a_s = a_s' = 0$ if $l \neq 1$, where
\begin{equation*}
\Psi_1 = \begin{pmatrix} 0 & z & 0 & 0 \\ 0 & 0 & 0 & 0 \\ 0 & x & 0 &
0 \\ 0 & y & 0 & 0 \end{pmatrix}, \Psi_1' = \begin{pmatrix} 0 & -w & -y
& x \\ 0 & 0 & 0 & 0 \\ 0 & 0 & 0 & 0 \\ 0 & 0 & 0 & 0 \end{pmatrix},
\end{equation*}
and
\begin{equation*}
P_2^0 = \begin{pmatrix} 0 & 0 & y & 0 \\ 0 & 0 & 0 & -1 \\ 0 & 0 & 0 &
0 \\ 0 & 0 & 0 & 0 \end{pmatrix}, P_3^0 = \begin{pmatrix} 0 & 0 & 0 & 0
\\ 0 & 0 & 0 & 0 \\ -1 & 0 & 0 & 0 \\ 0 & y & 0 & 0 \end{pmatrix}
\end{equation*}
The relation $2xy D_1 + z D_2 - w D_3 = 0$ gives the equation $2xy
P_1 + z P_2 - w P_3 = 0$ in $\enm_A(Y_l)$, i.e.
\begin{equation*}
\begin{pmatrix} 2(a_1-1)xy & bz-b'w & yz -b'y & b'x \\ 0 & 2(a_1-1) xy &
0 & -z \\ w & bx & 2a_1 xy & 0 \\ 0 & -yw + by & 0 & 2 a_1 xy
\end{pmatrix} = 0,
\end{equation*}
where $b = a_2z-a_3w$ and $b' = a_2'z-a_3'w$. By inspection, we see
that this is a contradiction for $1 \le l \le p$.

\emph{The module $B_1$ for $n$ even}. Let $P$ and $C$ be matrices
corresponding to graded endomorphisms of $L_0$ and $L_1$ of degree
$\omega$. Then $\deg p_{ij} = \deg e_j - \deg e_i + \omega$ and $\deg
c_{ij} = \deg f_j - \deg f_i + \omega$. In case $\omega = 0$, $P$ and
$C$ must be of the form
\begin{equation*}
P = \begin{pmatrix} p_{11} & 0 \\ 0 & p_{22} \end{pmatrix}, \; C =
\begin{pmatrix} c_{11} & 0 \\ 0 & c_{22} \end{pmatrix}
\end{equation*}
with $p_{ij}, c_{ij} \in k$ for $1 \le i,j \le 2$. In case $\omega =
n-3$, $P$ and $C$ must be of the form
\begin{align*}
P & = \begin{pmatrix} 0 & p_{12} y^{n-3} + p_{12}' x y^{p-1} \\ p_{21}
& 0 \end{pmatrix} \\ C & = \begin{pmatrix} 0 & c_{12} y^{n-3} + c_{12}'
x y^{p-1} \\ c_{21} & 0 \end{pmatrix}
\end{align*}
with $p_{ij}, p_{ij}', c_{ij}, c_{ij}' \in k$ for $1 \le i,j \le 2$.
For $s=1,4,5$, we must solve the equation $( D_s(\textsf B_1) ) =
\textsf B_1 C_s - P_s \textsf B_1$, where $D_s$ is the derivation of
$A$ mentioned in the proof of theorem \ref{t:simple-dim3}, and $P_s$
and $C_s$ have the above form (with $\omega = 0$ for $s = 1$ and
$\omega = n-3$ for $s = 4,5$). This gives $P_1 = P_1^0 + a_1 \Phi_0$
with $a_1 \in k$, where $\Phi_0 = I_2$ and
\begin{equation*}
P_1^0 = \begin{pmatrix} -1 & 0 \\ 0 & 0 \end{pmatrix}
\end{equation*}
Moreover, $P_s = P_s^0$ for $s = 4,5$, where $\gamma = (n-2) y^{n-3}$
and
\begin{equation*}
P_4^0 = \begin{pmatrix} 0 & -\gamma \\ 0 & 0 \end{pmatrix}, \; P_5^0 =
\begin{pmatrix} 0 & 0 \\ 1 & 0 \end{pmatrix}
\end{equation*}
Let $\beta = x^2 + (n-1) y^{n-2}$. The relation $\beta D_1 + zD_4 -
wD_5 = 0$ gives the equation $\beta P_1 + z P_4 - w P_5 = 0$ in
$\enm_A(B_1)$, i.e.
\begin{equation*}
\begin{pmatrix} - \beta + \delta & -z \gamma \\ -w & \delta
\end{pmatrix} = 0,
\end{equation*}
where $\delta = a_1 \beta$. By inspection, we see that this is a
contradiction.

\emph{The module $B_2$ for $n$ even}. We see that $B_2 \cong B_1^\vee$
using lemma \ref{l:mfid}. It follows from lemma \ref{l:dual} and the
computation above for $B_1$ that the module $B_2$ cannot admit a
$\g$-connection.

\emph{The module $C_-$ for $n$ even}. Let $P$ and $C$ be matrices
corresponding to graded endomorphisms of $L_0$ and $L_1$ of degree
$\omega$. Then $\deg p_{ij} = \deg e_j - \deg e_i + \omega$ and $\deg
c_{ij} = \deg f_j - \deg f_i + \omega$. In case $\omega = 0$, $P$ and
$C$ must be of the form
\begin{equation*}
P = \begin{pmatrix} p_{11} & 0 \\ 0 & p_{22} \end{pmatrix}, \; C =
\begin{pmatrix} c_{11} & 0 \\ 0 & c_{22} \end{pmatrix}
\end{equation*}
with $p_{ij}, c_{ij} \in k$ for $1 \le i,j \le 2$. In case $\omega =
n-3$, $P$ and $C$ must be of the form
\begin{align*}
P & = \begin{pmatrix} 0 & p_{12} y^{p-1} \\ p_{21} x + p_{21}' y^p & 0
\end{pmatrix} \\
C & = \begin{pmatrix} 0 & c_{12} y^{p-1} \\ c_{21} x + c_{21}' y^p & 0
\end{pmatrix}
\end{align*}
with $p_{ij}, p_{ij}', c_{ij}, c_{ij}' \in k$ for $1 \le i,j \le 2$.
For $s = 1,4,5$, we must solve the equation $( D_s(\textsf C_-) ) =
\textsf C_- C_s - P_s \textsf C_-$, where $D_s$ is the derivation of
$A$ mentioned in the proof of theorem \ref{t:simple-dim3}, and $P_s$
and $C_s$ have the above form (with $\omega = 0$ for $s = 1$ and
$\omega = n-3$ for $s = 4,5$). This gives $P_1 = P_1^0 + a_1 \Phi_0$
with $a_1 \in k$, where $\Phi_0 = I_2$ and
\begin{equation*}
P_1^0 = \begin{pmatrix} -1 & 0 \\ 0 & 0 \end{pmatrix}
\end{equation*}
Moreover, $P_s = P_s^0$ for $s = 4,5$, where
\begin{equation*}
P_4^0 = \begin{pmatrix} 0 & ip y^{p-1} \\ 0 & 0 \end{pmatrix}, \;
P_5^0 = \begin{pmatrix} 0 & 0 \\ -x+i(p+1)y^p & 0 \end{pmatrix}
\end{equation*}
Let $\beta = x^2 + (n-1) y^{n-2}$. The relation $\beta D_1 + zD_4 -
wD_5 = 0$ gives the equation $\beta P_1 + z P_4 - w P_5 = 0$ in
$\enm_A(C_-)$, i.e.
\begin{equation*}
\begin{pmatrix} (a_1-1) \beta & ip y^{p-1}z \\ xw-i(p+1)y^pw & a_1 \beta
\end{pmatrix} = 0
\end{equation*}
By inspection, we see that this is a contradiction.

\emph{The module $C_+$ for $n$ even}. Let $P$ and $C$ be matrices
corresponding to graded endomorphisms of $L_0$ and $L_1$ of degree
$\omega$. Then $\deg p_{ij} = \deg e_j - \deg e_i + \omega$ and
$\deg c_{ij} = \deg f_j - \deg f_i + \omega$. In case $\omega = 0$,
$P$ and $C$ must be of the form
\begin{equation*}
P = \begin{pmatrix} p_{11} & 0 \\ 0 & p_{22} \end{pmatrix}, \; C =
\begin{pmatrix} c_{11} & 0 \\ 0 & c_{22} \end{pmatrix}
\end{equation*}
with $p_{ij}, c_{ij} \in k$ for $1 \le i,j \le 2$. In case $\omega =
n-3$, $P$ and $C$ must be of the form
\begin{align*}
P & = \begin{pmatrix} 0 & p_{12} y^{p-1} \\ p_{21} x + p_{21}' y^p & 0
\end{pmatrix} \\
C & = \begin{pmatrix} 0 & c_{12} y^{p-1} \\ c_{21} x + c_{21}' y^p & 0
\end{pmatrix}
\end{align*}
with $p_{ij}, p_{ij}', c_{ij}, c_{ij}' \in k$ for $1 \le i,j \le 2$.
For $s = 1,4,5$, we must solve the equation $( D_s(\textsf C_+) ) =
\textsf C_+ C_s - P_s \textsf C_+$, where $D_s$ is the derivation of
$A$ mentioned in the proof of theorem \ref{t:simple-dim3}, and $P_s$
and $C_s$ have the above form (with $\omega = 0$ for $s = 1$ and
$\omega = n-3$ for $s = 4,5$). This gives $P_1 = P_1^0 + a_1 \Phi_0$
with $a_1 \in k$, where $\Phi_0 = I_2$ and
\begin{equation*}
P_1^0 = \begin{pmatrix} -1 & 0 \\ 0 & 0 \end{pmatrix}
\end{equation*}
Moreover, $P_s = P_s^0$ for $s = 4,5$, where
\begin{equation*}
P_4^0 = \begin{pmatrix} 0 & -ip y^{p-1} \\ 0 & 0 \end{pmatrix}, \;
P_5^0 = \begin{pmatrix} 0 & 0 \\ -x-i(p+1)y^p & 0 \end{pmatrix}
\end{equation*}
Let $\beta = x^2 + (n-1) y^{n-2}$. The relation $\beta D_1 + zD_4 -
wD_5 = 0$ gives the equation $\beta P_1 + z P_4 - w P_5 = 0$ in
$\enm_A(C_+)$, i.e.
\begin{equation*}
\begin{pmatrix} (a_1-1) \beta & -ip y^{p-1}z \\ xw + i(p+1)y^p w & a_1
\beta \end{pmatrix} = 0
\end{equation*}
By inspection, we see that this is a contradiction.

\emph{The modules $D_-$ and $D_+$ for n even}. We see that $D_- \cong
C_+^\vee$ and that $D_+ \cong C_-^\vee$ using lemma \ref{l:mfid}. It
follows from lemma \ref{l:dual} and the computations above for $C_-$
and $C_+$ that the modules $D_-$ and $D_+$ cannot admit
$\g$-connections.

\bibliographystyle{amsplain}
\bibliography{eeriksen}

\end{document}